\documentstyle[12pt,leqno,amstex,amssymb,amsfonts,eucal]{article}
\begin{document}
\oddsidemargin= 10mm
\topmargin= 35mm
\textwidth=170mm
\textheight=270mm
\pagestyle{plain}
\newcounter{r}
\newcommand{\Ker}{Ker}
\newcommand{\Gal}{Gal}
\newcommand{\Mat}{Mat}
\newcommand{\sign}{sign}
\newcommand {\Log}{Log}
\newcommand{\Res}{Res}
\newcommand{\Tr}{Tr}
\newcommand{\Nm}{Nm}
\begin{center}
{\large\bf  On the Diophantine Approximations of logarithms
in cyclotomic fields.}\end{center}
\begin{center}
\vskip.4pt
{\large\bf L.A. Gutnik}
\end{center}
\vskip10pt
\small Abstract. 
I present here the proofs of results, which are obtained in my papers
" On the linear forms, whose
coefficients are linear combinations with algebraic coefficients
of logarithms of algebraic numbers,"
 VINITI, 1996,  1617-B96, pp. 1 -- 23 (in Russian),
and
" On the systems of linear forms, whose
coefficients are linear combinations with algebraic coefficients
of logarithms of algebraic numbers,"
 VINITI, 1996,  2663-B96, pp. 1 -- 18 (in Russian),
\vskip 10pt
\normalsize
\vskip.10pt

 \rm Let $T$ is a real number, $\Delta,m$ and $n,$ are positive integers,
 $\Delta$ is greater or equal than $2,$
$K_m={{\mathbb Q}[\exp(2\pi i/m)]}$
is a cyclotomic field, ${\mathbb Z}_{K_m}$ is the ring of all the integers
of the field $K_m,\,\Lambda(n)$ is the Mangold's function,
 $\epsilon^2=\epsilon.$
Let $\Lambda_0(m)=0,$ if $m$ is odd and $\Lambda_0(m)=\Lambda(m/2),$
if $m$ is even. Let further $\omega_1(m)=(m-1)/2,$ if $m$ is odd,
$\omega_1(m)=m/2-2,$ if $m\equiv2(\mod4)$ and $\omega_1(m)=m/2-1,$
if $m\equiv0(\mod4).$
Let
\begin{equation}
w_\Delta(T)=\sqrt{\frac{\sqrt{(\Delta^2(3-T^2)+1)^2+16\Delta^4T^2}+
\Delta^2(3-T^2)+1}2},
\label{eq:0}\end{equation}
and the values $V_{\Delta}^\ast,\,V_{\Delta}(m),\,l_\Delta(\epsilon,T),\,
g_{\Delta,k}(m)$ and $h_{\Delta,k}(m)$
are defined by the equalities
\begin{equation}
V_\Delta^\ast=(\Delta+1)+
\log((\Delta-1)^{(\Delta-1)/2}(\Delta+1)^{(\Delta+1)/2}\Delta^{-\Delta})\,+
\label{eq:1}\end{equation}
$$\frac\pi2\sum\limits_{\mu=0}^1(1-2\mu)
\sum\limits_{\kappa=1}^{[(d-1)/2]+\mu}
\cot\left(\frac{\pi\kappa}{d-1+2\mu}\right),
$$
\begin{equation}
V_\Delta(m)=V^\ast+(\Delta+1)\Lambda_0(m)/\phi(m),
\label{eq:3}\end{equation}
\begin{equation}
l_\Delta(\epsilon,T)=
-log\left(4(\Delta+1)^{\Delta+1}(1-1/\Delta)^(\Delta-1)\right)\,+
\label{eq:2}\end{equation}
$$\frac12\log\left(\left(2\Delta+(-1)^\epsilon w_\Delta(T)+(\Delta+1)\right)^2+
T^2\Delta^2\left(1+\frac{(-1)^\epsilon2\Delta}{w_\Delta(T)}\right)^2\right)+$$
$$\frac12
\log\left(\left(2\Delta+(-1)^\epsilon w_\Delta(T)-(\Delta+1)\right)^2+
T^2\Delta^2\left(1+\frac{(-1)^\epsilon2\Delta}{w_\Delta(T)}\right)^2\right)+$$
\newpage
\pagestyle{headings}
\topmargin= -15mm
\textheight=250mm
\markright{\footnotesize\bf L.A.Gutnik,
 On the Diophantine Approximations of logarithms in cyclotomic fields.}
$$\frac{(\Delta-1)}2
\log\left(\left(2\Delta+(-1)^\epsilon w_\Delta(T)\right)^2+
T^2\Delta^2\left(1+\frac{(-1)^\epsilon2\Delta}{w_\Delta(T)}\right)^2\right),$$
\begin{equation}
g_{\Delta,\epsilon}(m)=
(-1)^\varepsilon(l_\Delta(\epsilon,\tan(\pi\omega_1(m)/m)+V_\Delta(m))),
\label{eq:2a0}\end{equation}
\begin{equation}
h_{\Delta}(m)=-V_\Delta(m)-l_\Delta(1,\tan(\pi/m)),
\label{eq:2a}\end{equation}
where $m\ne 2,\,k=0,1.$ Let, finally, the values $\beta(d,m)$
and $\alpha(\Delta,m)$ are defined by the equalities
$$\beta(\Delta,m)=g_{d,0}(m)/h_{\Delta}(m),\,
\alpha(\Delta,m)=\beta(\Delta,m)-1+g_{\Delta,1}(m)/h_{\Delta}(m).$$
\vskip1pt
\bf Theorem. \it Let $m$ is a positive integer different from one, two
 and six, $$\Delta\in\{5,7\}.$$
 Then
\begin{equation}
h_\Delta(m)>0
\label{eq:4}\end{equation}
 and for each $\varepsilon>0$
 there exists $C_{\Delta,m}(\varepsilon)>0$ such that
\begin{equation}
\max_{\sigma\in Gal(K/{\mathbb Q})}
(\vert q^\sigma\log((2+\exp(2\pi i/m))^\sigma)-p^\sigma\vert)\ge
\label{eq:5}\end{equation}
$$C_{\Delta,m}(\varepsilon)(\max_{\sigma\in\Gal(K_m/{\mathbb Q})}
(\vert q^\sigma\vert)^{-\alpha(\Delta,m)-\varepsilon},$$
where $p\in{\mathbb Z}_{K_m}$ and $q\in{\mathbb Z}_{K_m}\diagdown \{0_{K_m}\};$
 moreover, for any $q\in{\mathbb Z}_{K_m}\diagdown \{0_{K_m}\}$
and any $\varepsilon>0$ there exists $C^\ast_{\Delta,m}(q,\varepsilon)>0$
 such that
\begin{equation}
b^{{\beta(\Delta,m)+\varepsilon}}\max_{\sigma\in Gal(K/{\mathbb Q})}
(\vert q^\sigma b\log((2+\exp(2\pi i/m))^\sigma)-p^\sigma\vert)\ge
\label{eq:6}\end{equation}
$$C^\ast_{\Delta,m}(q,\varepsilon),$$
where $p\in{\mathbb Z}_{K_m},\,b\in{\mathbb N}.$

 For the proof I use the same method, as in [\ref{r:ag}] -- [\ref{r:ci2}].
I must work on the Riemann surface $\mathfrak F$  of the function
$\Log (z)$ and identify it with the direct product of the multiplicative group
${\mathbb R}_+^\ast= \{r\in{\mathbb R}\colon r>0\}$ of all the positive
 real numbers
with the operation $\times,$ not to be written down explicitly as usual,
 and the additive group ${\mathbb R}$ of all the real numbers, so that
$$z_1z_2=(r_1r_2,\phi_1+\phi_2)$$ for any two points
$z_1=(r_1,\phi_1)$ and $z_2 = (r_2,\phi_2)$ on $\mathfrak F.$
I will illustrate the appearing situations on the half plain
$(\phi,r),$ where $r>0.$

 For each $z=(r,\phi)\in{\mathfrak F},$  let
$$\theta_0 (z)=r\exp {i\phi},\,\Log(z)=\ln(r)+i\phi,\,\eta^\ast_\alpha(z)
= (r,\phi-\alpha),$$
\newpage
where $\alpha\in{\mathbb R}.$
Clearly,  $\Log(z_1z_2)=Log(z_1)+\Log(z_2)$ for any $z_1\in {\mathfrak F}$
$z_2\in{\mathfrak F}.$
Let $\rho(z_1,z_2)=\vert\Log(z_1)-\Log(z_2)\vert,$
where $z_1\in{\mathfrak F}$ and $z_2\in {\mathfrak F};$ clearly,
 $({\mathfrak F},\rho)$ is a metric space.
 Clearly,
$\rho(zz_1,zz_2)=\rho(z_1,z_2)$ for any $z_1,z_2$ and $z$ in ${\mathfrak F}.$
Clearly, $\theta_0(z)=\exp(\Log(z))$ for any $z\in{\mathfrak F}.$
Clearly, for any $\alpha\in{\mathbb R}$ the map $z\to\eta^\ast_\alpha(z)$
is the bijection of ${\mathfrak F}$ onto ${\mathfrak F}$ and
$$\theta_0((\eta^\ast_\alpha)^m(z))=\exp(-im\alpha)\theta_0(z)$$
for each $z=(r,\phi)\in{\mathfrak F},\,\alpha\in{\mathbb R}$
 and $m\in{\mathbb Z}.$ Clearly, the group ${\mathfrak F}$ may be considered
 as ${\mathbb C}$-linear space, if for
 any $z\in{\mathfrak F}$ and any $s\in{\mathbb C}$ we let
$$
z^s=(\vert\exp(s\Log(z))\vert,\Im(s\Log(z)).
$$
 Let us fix a domain $D$ in
${\mathfrak F}.$ Let $f(z)=f^\wedge (r,\phi)$ for a complex-valued
 function $f(z)$ on $D,$ It is well known that $f(z)$ is holomorphic in $D$
 if the complex-valued function $f^\wedge (r,\phi)$ of two real variables $r$
 and $\phi$ has continuous partial derivatives in $D,$ and
 the Cauchy-Riemann conditions
\begin{equation}
r(((\partial/\partial r)f^\wedge)(r,\phi))=
-i((\partial/\partial\phi)
f^\wedge)(r,\phi)):=\label{eq:a}\end{equation}
$$(\delta f)(z):=
 \theta_0(z)((\partial/\partial z)f)(z))
$$
are satisfied  for every point $z=(r,\phi)\in D$. The equalities (\ref{eq:a})
 determine a differentiations $\frac{\partial}{\partial z}$
 and $\delta=\theta_0(z)\frac{\partial}{\partial z}$ on the ring of all
 the holomorphic in the domain D functions. In particular,
 the function $\Log(z)$ is holomorphic on $\mathfrak F$ and we have
 the equalities
$$\left((\partial/\partial z)\Log\right)(z)=\theta_0(z^{-1}),\,
(\delta\Log)(z)=1.$$
For the proof I use the functions of C.S.Mejer.
Let $\Delta\in{\mathbb N}+1,\;\delta_0=1/\Delta,$
$$
\gamma_1=(1-\delta_0)/(1+\delta_0),\quad
d_{l}=\Delta+(-1)^l, \quad l = 1, 2.
$$
To introduce the first of my auxiliary function $f_1(z,\nu),$ I use
 the auxiliary set
$$
\Omega_0=\{z\in{\mathfrak F}\,\colon\,\vert z\vert\le1\}.
$$
I  prove that, for each $\nu \in {\mathbb N},$ the function
$f_1(z,\nu)$ belongs to the ring ${\mathbb Q}[\theta_0(z)];$ therefore using
the principle of analytic continuation we may regard it as being
defined in ${\mathfrak F}.$ For $\nu\in{\mathbb N},$ let
\begin{equation}
f_1(z,\nu) =
 -(-1)^{\nu(\Delta+1)}G_{2, 2}^{(1, 1)}
\left(z\bigg\vert\matrix -\nu d_1,
 &\!\! 1+\nu d\sb2\\ 0, &\!\! \nu\\\endmatrix\right)
\label{eq:b}\end{equation}
$$= -(-1)^{\nu(\Delta+1)}\frac1{2\pi i}\int\limits_{L_1}g_{2, 2}^{(1,1)}(s)ds,
$$
where
$$
g_{2,2}^{(1,1)}(s)= \theta_0(z^s)\Gamma(-s)
\Gamma(1+d_1\nu+s)/(\Gamma (1-\nu +s)\Gamma
(1+d_2\nu-s))
$$
and the curve $L_1$ passes from $+\infty$ to $+\infty$ encircling
the set ${\mathbb N}-1 $ in the negative direction, but not including any
point of the set $-{\mathbb N}.$ So, for the parameters of the
Meyer's functions we have
$$p=q=2,\,m=n= 1,\,a_1=-\nu d_1,\,a_2=1+\nu d_2,\,b_1=0,\,b_2=\nu,$$
$$\Delta^\ast=\left(\sum\limits_{k=1}^qb_k\right)-
\sum\limits_{j=1}^pa_j=-\nu-1<-1,$$
and, since we take $\vert z \vert \le  1,$ convergence conditions
of the integral in (\ref{eq:b}) hold. To compute the function $f_1(z,\nu),$
 we use the following formula
\begin{equation}
G = (-1)^k\sum\limits_{s\in S_k}\Res(g;s),
\label{eq:c}\end{equation}
where $k=1,\,G$ denotes the integral (\ref{eq:b}) with $L = L\sb k,$
$g$ denotes the integrand of the integral (\ref{eq:b}),
 $S_k$ denotes
the set of all the unremovable singularities of $g$ encircled by $L_k,$
and
$\Res(g;s)$ denotes the residue of the function $g$ at the point $s.$
Then we obtain the equlity
$$
f_1(z,\,\nu) =
$$
$$ (\nu d\sb 1)!/(\nu \Delta)!
z^\nu(-1)^{\nu\Delta}\sum\limits_{k=0}^{\nu\Delta}(-\theta_0(z))^k
 \binom {\nu \Delta  }k\binom {\nu \Delta + k  }{\nu d \sb 1}.$$

Therefore, as it has been already remarked, using the principle of
analytic continuation we may regard it as being defined in
$\mathfrak F.$
Let
$$
\Omega_1=\{z\in{\mathfrak F}\,\colon\,\vert z\vert\ge 1\}.
$$
Now, let me introduce my second auxiliary function defined
for $z\in\Omega_1.$
For $\nu\in{\mathbb N},$ let
\begin{equation}
f_2(z,\nu) =
 -(-1)^{\nu(\Delta+1)}G_{2, 2}^{(2, 1)}
\left(z\bigg\vert\matrix -\nu d_1,
 &\!\! 1+\nu d\sb2\\ 0, &\!\! \nu\\
\endmatrix\right)=
\label{eq:e}\end{equation}
$$
-(-1)^{\nu(\Delta+1)}\frac1{2\pi i}\int\limits_{L_2}g_{2, 2}^{(2,1)}(s)ds,
$$
where
$$
g_{2,2}^{(2,1)}(s)= \theta_0((\eta_\pi(z))^s)\Gamma(-s)\Gamma (\nu-s)
\Gamma(1+d_1\nu+s)/\Gamma (1+d_2\nu-s).
$$
and the curve $L_2$ passes from $-\infty$ to $-\infty$ encircling
the set $-{\mathbb N}$ in the positive direction,
 but not including any point of the set ${\mathbb N}-1.$
  So, for the parameters of the
Meyer's functions we have $$p=q=m=2,\,n= 1,\,
a_1=-\nu d_1,\,a_2=1+\nu d_2, \,
b_1=0,\,b_2=\nu,$$
$$\Delta^\ast=\left(\sum\limits_{k=1}^qb_k\right)-
\sum\limits_{j=1}^pa_j=-nu-1<-1,$$
and, since we take $\vert z \vert \ge  1,$ convergence conditions
of the integral in (\ref{eq:e}) hold. To compute the function $f_2(z,\nu),$
 we use the formula (\ref{eq:c})
where $k=2,\,G$ denotes the integral in (\ref{eq:e}) with $L = L\sb k,$
$g$ denotes the integrand of the integral in (\ref{eq:e}),
 $S_k$ denotes
the set of all the unremovable singularities of $g$ encircled by $L_k,$
and
$\Res(g;s)$ denotes the residue of the function $g$ at the point $s.$
Then we obtain the equality
\begin{equation}
f_2 (z, \nu )(\nu \Delta )!/(\nu d\sb1)! =
(-1)^\nu \!\!\sum \limits \sb {t=\nu + 1} \sp \infty R_0(t;\nu)
\theta_0(z^{-t+\nu}),
\label{eq:f}\end{equation}
where
$$
 R_0(t;\nu)=(\nu \Delta )!/(\nu d\sb1)!
\left(\prod\limits_{\kappa=\nu+1}^{\nu\Delta}(t-\kappa)\right)
\prod\limits_{\kappa=0}^{\nu\Delta}(t+\kappa)^{-1}.
$$
Let further
\begin{equation}\label{eq:h}
f_k^\ast(z,\,\nu )=f_k (z,\,\nu )(\nu \Delta )!/(\nu d\sb1)!,
\end{equation}
where $k=1,\,2.$
Expanding the function $R_0(t;\nu)$ into partial fractions,
we obtain the equality
$$
R_0(t;\nu)=\sum\limits_{k=0}^{\nu\Delta}
\alpha^\ast_{\nu, k}/(t + k)
$$
with
\begin{equation}
\alpha\sp\ast\sb{\nu, k} = (-1)\sp {\nu + \nu \Delta + k}
\binom{\nu \Delta} {k}\binom{\nu \Delta+k}
{\nu \Delta  - \nu},
\label{eq:i}\end{equation}
where $k=0,\,\ldots,\,\nu\Delta.$
It follows from (\ref{eq:e}), (\ref{eq:f}), (\ref{eq:h}) and (\ref{eq:i})
that
\begin{equation}
f^\ast_2(z,\nu)=(-\theta_0(z))^\nu
\sum\limits_{t=1+\nu}^{+\infty}(\theta_0(z))^{-t}R\sb 0(t;\nu)=
\label{eq:aa}\end{equation}
$$
=(-\theta_0(z))^\nu\sum \limits_{t=1+\nu}^{+\infty}(\theta_0(z))^{-t-k+k}
\sum\limits_{k=0}^{\nu\Delta}\alpha^\ast_{\nu,k}/(t+ k)
$$
$$
=(-\theta_0(z))^\nu\sum \limits_{t=1+\nu}^{+\infty}
((\theta_0(z))^{-t-k}/(t+k))
\sum\limits_{k=0}^{\nu\Delta}\alpha^\ast_{\nu,k}(\theta_0(z))^{k}=
$$
$$
(-\theta_0(z))^\nu
\sum\limits_{k=0}^{\nu\Delta}\alpha^\ast_{\nu,k}(\theta_0(z))^{k}
\sum \limits_{\tau=1+\nu+k}^{+\infty}
((\theta_0(z))^{-\tau}/\tau))=
$$
$$
=\alpha^\ast (z;\nu)(-\log(1-1/\theta_0(z))) - \phi^\ast(z;\nu),
$$
where $\log(\zeta)$ is a branch of $\Log(\zeta)$
 with $\vert\arg(\zeta)\vert<\pi,$
\begin{equation}
\alpha^\ast(z;\nu)=(-(\theta_0(z))^\nu\sum\limits_{k=0}^{\nu\Delta}
\alpha^\ast_{\nu,\,k}(\theta_0(z))^k=f^\ast_1 (z;\nu),
\label{eq:ab}\end{equation}
\begin{equation}
\phi^\ast(z;\nu)=(-\theta_0(z))^\nu
\sum\limits_{k=0}^{\nu\Delta}\alpha^\ast_{\nu,k}(\theta_0(z))^{k}
\sum \limits_{\tau=1}^{\nu+k}
((\theta_0(z))^{-\tau}/\tau))=
\label{eq:ac}\end{equation}
$$
(-\theta_0(z))^\nu\sum \limits_{\tau=1}^\nu((\theta_0(z))^{-\tau}
\alpha^\ast(z;\nu)/\tau+
$$
$$
(-\theta_0(z))^\nu
\sum\limits_{k=0}^{\nu\Delta}\alpha^\ast_{\nu,k}(\theta_0(z))^{k}
\sum \limits_{\tau=1+\nu}^{\nu+k}
((\theta_0(z))^{-\tau}/\tau)).
$$
The change of order of summation by passage to (\ref{eq:aa})
 is possible, because the series in the second sum in (\ref{eq:aa})
is convergent, if $\vert z\vert\ge1$ and $\theta_0(z)\ne1.$
Since
$$\deg_t\left(\prod\limits_{\kappa=\nu+1}^{\nu\Delta}(t-\kappa)\right)-
\deg_t\left(\prod\limits_{\kappa=0}^{\nu\Delta}(t+\kappa)\right)=-\nu-1,
$$
it follows that
$$
\alpha^\ast(1;\nu)=\Res(R_0(t;\nu);t=\infty)=0
$$
So in the domain $D_0 = \{z\in{\mathfrak F}\colon \vert z\vert>1$
the funcion $f^\ast_2(z,\nu)$ coincides with the functin
\begin{equation}
f^\ast_0(z,\nu)=\alpha^\ast (z;\nu)(-\log(1-1/\theta_0(z)))
 - \phi^\ast(z;\nu),
\label{eq:aa1}\end{equation}
The form (\ref{eq:aa1}) may be used for various
applikations. Espeshially it is pleasant, when both $1/\theta_0(z)$
and $\alpha^\ast(z;\nu)$ for some $z$ is integer algebraic number.
The following Lemma corresponds to this remark.

\bf Lemma 1. \it Let $m\in{\mathbb N},\,m>2$
$m\ne2p^\alpha,$ where $p$ run over the all the prime numbers and $\alpha$
run over $\mathbb N.$ Then $1+\exp(2\pi i/m)$ belongs to the group of
 the units of the field $K_m.$  If $m=2p^\alpha,$
where $p$ is a prime number and $\alpha\in\mathbb N,$ then the ideal
${\mathfrak l}=(1+exp(2\pi i/m))$   is a prime ideal
 in the field $K_m,$ and ${\mathfrak l}^{\phi(m)}=(p).$

\bf Proof. \rm Let polynomial $\Phi_m(z)$ is irreducible over ${\mathbb Q},$
has the leading coefficient equal to one and $\Phi_m(\exp(2\pi i/m))=0.$
Let $\Lambda(n),$ as usual, denotes the Mangold's function.
Since (see, for example, [\ref{r:cb1}], end of the chapter 3)
$$\Phi_m(z)=\prod\limits_{d\vert m}(z^{m/d}-1)^{\mu(d)},$$
it follows that
$$\Phi_m(-1)=(-2)^{\left(\sum\limits_{d\vert m}\mu(d)\right)}=1,$$
if $m\in 1+2{\mathbb N},$
$$\Phi_m(z)=\prod\limits_{d\vert (m/2)}
(((z)^{m/(2d)}-1)^{\mu(2)}((-z)^{m/d}-1)/((-z)-1))^{\mu(d)},$$
$$\Phi_m(-1)=\lim\limits_{z\to-1}\prod\limits_{d\vert (m/2)}
(((-z)^{m/d}-1)/((-z)^-1))^{\mu(d)}\times$$
$$(-2)^{\mu(2)\left(\sum\limits_{d\vert (m/20}\mu(d)\right)}=$$
$$\exp\left(\sum\limits_{d\vert (m/2)}\ln(m/(2d)){\mu(2d)}\right)=
\exp(\Lambda(m/2)),$$
if $m\in 2(1+2{\mathbb N}),$
$$\Phi_m(z)=\prod\limits_{d\vert (m/2)}
(((-z)^{m/d}-1)/((-z)-1))^{\mu(d)},$$
and
$$\Phi_m(-1)=\lim\limits_{z\to-1}\prod\limits_{d\vert(m/2)}
(((-z)^{m/d}-1)/((-z)-1))^{\mu(d)}=$$
$$\exp\left(\sum\limits_{d\vert m/2}\ln(m/(2d)){\mu(d)}\right)
=\exp(\Lambda(m/2)),$$
if $m\in 4{\mathbb N}.$
If $m=2p^\alpha$ with $\alpha\in{\mathbb N},$
 then $\Phi_m(-1)=\exp(\Lambda(m/2))=p,$
and ideals ${\mathfrak l}_k=(1+exp(2\pi ik/m))$, where $(k,m)=1,$
 divide each other and in the standard equality $efg=n$
(see,  [\ref{r:cb1}], chapter 3, section 10)
 we have
$$e=n=\phi(m),\ f=g=1.\,\blacksquare.$$
In connection with the above remark and with the Lemma 1,
the following case is interesting for us:
\begin{equation}
\theta_0(z)=(-\rho)(1+\exp(-i\beta))=
-(\rho\exp(i\beta/2))/(2cos(\beta/2))=
\label{eq:aa3}\end{equation}
$$-(\rho\exp(i\psi))(2cos(\psi))=
-(1+i\tan(\psi))/2$$
with $\rho>2/3,\, \vert\beta\vert<\pi$ and $-\pi/2<\psi=\beta/2<\pi/2;$
 then
$$
\Re(1-1/\theta_0(z))=\Re(2+\exp(i\beta)/\rho)>1/2,
$$
 and we have no problems with $\log(1-1/\theta_0(z)).$
Of course, according to the Lemma 1, the case $\rho=1$  is interesting
 especially. So, we will take further
\begin{equation}
z=(\rho/(2\cos(\psi)),\psi-\pi)=
(\rho/(-2cos(\theta),\,\theta),
\label{eq:aa4}\end{equation}
where $\rho>2/3,\,\vert\psi\vert<\pi/2$ and $-3\pi/2<\theta=\psi-\pi<-\pi/2;$
 clearly, the function (\ref{eq:aa1}) is analytic in the domain
$$
D_1=
\{z=(\rho(2\cos(\psi))^{-1},\psi-\pi))\colon\rho>2/3,\,-\pi/2<\psi<\pi/2\}=
$$
$$
\{z=((-2\rho\cos(\theta))^{-1},\theta))\colon\rho>2/3,\,
-3\pi/2<\theta<-\pi/2\}.
$$
Let
\begin{equation}
D_2(\delta_0) = \{z\in{\mathfrak F}\colon \vert z\vert>1+\delta_0/2\},\,
D_3=D_2(\delta_0)\cup D_1.
\label{eq:aa6}\end{equation}
So, the funcion $f^\ast_2(z,\nu)$ coincides with the function (\ref{eq:aa1})
in $D_2(\delta_0)\subset D_0.$
Since $D_2(\delta_0)\cap D_1\ne\emptyset,$ it follows that the join
 $D_3=D_2(\delta_0)\cup D_1$
of the domains $D_2(\delta_0)$ and $D_1$ is a domain in $\mathfrak F$ and
the function (\ref{eq:aa1}) is analytic in this domain.

The conditions, which imply the equality
\begin{equation}
(-1)^{m+p-n}\exp(-i\alpha)\theta_0(z)\times
\label{eq:ae}\end{equation}
$$\left(\left(\prod\limits_{j=1}^p(\delta+1-a_j)\right)
(G\circ\eta^\ast_\alpha)\right)(z)=
\left(\left(\prod\limits_{k=1}^q(\delta-b_k)\right)
(G\circ\eta^\ast_\alpha)\right)(z)$$
hold in our case for the Mejer's function
$$
G=G \sb {p, q} \sp{(m, n)}\left(z\bigg\vert\matrix a_1,&\ldots,&a_p\\
b_1,& \ldots,&b_q\\
\endmatrix\right).$$
 We have $p=q=2,\,m=n=1,\,\alpha=0$ for the function $f_1(z,\nu)$
 and the equation (\ref{eq:ae}) takes the form
$$
\theta_0(z)
((\delta+1+d_1\nu)(\delta-d_2\nu)f_1)(z,\nu)=
(\delta(\delta-\nu)f_1)(z,\nu)
$$
We have $p=q=m=2,\,n=1,\,\alpha=\pi$ for the function $f_2(z,\nu)$
 and the equation (\ref{eq:ae}) takes the form
$$
\theta_0(z)
((\delta+1+d_1\nu)(\delta-d_2\nu)f_2)(z,\nu)=
(\delta(\delta-\nu)f_2)(z,\nu).
$$
We see that both the  functions $f(z,\nu)=f^\ast_k(z,\nu),$ where $k=1,2$
satisfy to the same differential equation
\begin{equation}
\theta_0(z)
(\delta+1+d_1\nu)(\delta-d_2\nu)f(z,\nu)=
(\delta(\delta-\nu)f)(z,\nu).
\label{eq:ah}\end{equation}
in the domain $D_0.$
According to the general properties of the Mejer's functions
 we have the equality
\begin{equation}
\left(\prod\limits_{\kappa=1}^{\Delta-1}(\nu(\Delta-1)+\kappa)\right)
\prod\limits_{\kappa=1}^{d_2}(\delta-d_2\nu-\kappa)
f^\ast_k(z,\,\nu+1) =
\label{eq:ai}\end{equation}
$$
\left(\prod\limits_{\kappa=1}^{\Delta}(\nu\Delta+\kappa)\right)(\delta-\nu)
\prod\limits_{\kappa=1}^{d_1}(\delta+d_1\nu+\kappa)f^\ast_k(z,\,\nu),
$$
where $k=1,2$ and $z\in D_0.$
Since $f^\ast_0(z,\nu)$ and polynomial $f^\ast_1(z,\nu)$ are analytic
in the domain $D_0\cup D_1,$ and $f^\ast_0(z,\nu)$
coincides with $f^\ast_2(z,\nu),$ it follows that the equations
(\ref{eq:ah}) and (\ref{eq:ai}) hold in $D_0\cup D_1$ for $k=0,1.$

Let
\begin{equation}
D^\vee (w,\eta)=(\eta+1)(\eta+\gamma_1)-
2(1+\gamma_1)w\eta,
\label{eq:bj}\end{equation}
\begin{equation}
D^\wedge (z,\eta)=D^\vee (\theta_0(z),\eta),
\label{eq:bj1}\end{equation}
where, in view of (\ref{eq:aa3}),
\begin{equation}
w=\theta_0(z)=-r\exp(i\psi),\,r=1/(2cos(\psi)),\,\vert\psi\vert<\pi/2.
\label{eq:bj2}\end{equation}
In view of (\ref{eq:bj2}), the polynomial (\ref{eq:bj}) coincides
with the polynomial (1) in [\ref{r:cc}]. Let 
\begin{equation}
h^\sim(\eta)=
(\eta-1)(1-\delta_0)^{-d_1}(\eta+1)2^{-2} \
\eta^{d_1}.
\label{eq:ba}\end{equation}
 As in [\ref{r:bg}], we consider $\nu^{-1}$ as an independent variable
taking its values in the field ${\mathbb C}$ including $0.$
Let $F$ be a bounded closed subset of ${\mathfrak F} $ (in
particular, this compact $F$ may be an one-point set).
Let ${\mathfrak H}_0 (F)$ be the subring
 of all those functions in ${\mathbb Q}(w),$ which are well defined for every
$w\in\theta_0(F).$ For $\varepsilon\in(0,1)$, let
${\mathfrak H} (F,\varepsilon)$ be the subring of all those functions in
${\mathbb Q}(w,\nu^{-1})$, which are well defined for every
$(w, \nu^{-1})$ with
$w\in \theta_0(F),\,\vert\nu^{-1}\vert\le\varepsilon_0.$

\bf Lemma 2. \it Let $F$ be a closed bounded
subset of $D_0\cup D_1$ (in particular, $F$ may be an one-point set).
 Let further for any $ z\in F$ the polynomial (\ref{eq:bj1})
 has only simple roots and on the set of all the roots $\eta$ of
 the polynomial $D^\wedge (z,\eta)$ the map
\begin{equation}
 \eta\to h^\sim(\eta)
\label{eq:bb}\end{equation}
is injective. Then there is
$\varepsilon\in (0,1)$ such that, for
any $z\in F,\nu\in {\mathbb N}+[1\diagup\varepsilon]$,
the functions
$f^\ast_0(z,\nu),\,f^\ast_1(z,\nu)=\alpha^\ast (z;\nu)$ and
$\phi^\ast(z;\nu)$ are solutions of the difference equation
\begin{equation}
x(z,\nu+2)+\sum\limits_{j=0}^1 q^\ast_j(z,\nu^{-1})x(z,\nu+j) = 0,
\label{eq:bc}\end{equation}
moreover,
\begin{equation}
q^\ast_j (z,\,\nu^{-1})\in{\mathfrak H}(F,\,\varepsilon)
\label{eq:bd}\end{equation}
for $j = 0,\,1,$ and trinomial
\begin{equation}
 w^2+\sum\limits_{j=0}^1q^\ast_j (z,0)w^j
\label{eq:be}\end{equation}
coincides with
\begin{equation}
\prod\limits_{k=0}^1(w-h(\eta_k)),
\label{eq:bd1}\end{equation}
if
$$\prod\limits_{k=0}^1(w-\eta_k),$$
coincides with $D^\vee (w,\eta)$ from (\ref{eq:bj}).

\bf Proof. \rm Proof may be found in [\ref{r:bg}]. $\blacksquare$

This Lemma shows the importance of the properties of the roots
 of the polynomial (\ref{eq:bj}).
In correspondence with (\ref{eq:aa4}) and with notations in [\ref{r:cc}], let
\begin{equation}
\rho>2/3,\,r=\rho/(2\cos(\psi)),\,t=\cos(\psi),\,
\vert\psi\vert<\pi/2.
\label{eq:be1}\end{equation}
Let
$
u=r^2,\delta_0\le1/2<2/3<\rho.
$
 Then
\begin{equation}
2\delta_0\le2/5<2/3<\rho<2\sqrt{u}=2r.
\label{eq:be3}\end{equation}
 Clearly,
$$(\partial/\partial\psi)r=(\rho\sin(\psi))/(2\cos^2(\psi))=
-2\rho(\sin(\psi)-1)-2\rho/(\sin(\psi)+1),$$
$$(\partial/\partial\psi)^2r=
(2\rho\cos(\psi))/(\sin(\psi)-1)^2
+(2\rho\cos(\psi))/(\sin(\psi)+1)^2>0,$$
if $\vert\psi\vert<\pi/2$
 In view of (3.1.10) in [\ref{r:bh}],
\begin{equation}\label{eq:bf}
\vert D_0(r,\psi,\delta_0 )\vert^2 =
r^4+r^2+(\delta_0/2)^4+\end{equation}
$$2r^2(\delta_0/2)^2(2t^2-1)+2r(r^2+(\delta_0/2)^2)t=$$
$$u^2+u+(\delta_0/2)^4+(\delta_0/2)^2(\rho^2-2u)+\rho(u+(\delta_0/2)^2)=$$
$$u^2+u(\rho+1-(\delta_0)^2/2)+
(\delta_0/2)^2(\rho^2+\rho+(\delta_0/2)^2),$$
\begin{equation}
\vert R_0(r,\psi,\delta_0)\vert^2=\vert D_0(r,\psi,\delta_0)\vert=
\label{eq:bg}\end{equation}
$$\sqrt{u^2+u(\rho+1-(\delta_0)^2/2)+
(\delta_0/2)^2(\rho^2+\rho+(\delta_0/2)^2)}.$$
 In view of (3.1.41) - (3.1.43) in [\ref{r:bh}] and (\ref{eq:bg}),
\begin{equation}
p_1 =8(\vert R^\ast_0(r,\psi,\delta_0)\vert^2+
\vert R_0(r,\psi,\delta_0)\vert^2)/(1+\delta_0)^2=
\label{eq:bh}\end{equation}
$$
8(r^2+rt+1/4+\vert D_0(r,\psi,\delta_0)\vert)/
(1+\delta_0)^2=8(1+\delta_0)^{-2}\times
$$
$$
\left(u+\rho/2+1/4+
\sqrt{u^2+u(\rho+1-(\delta_0)^2/2)+
(\delta_0/2)^2
(\rho^2+\rho+(\delta_0/2)^2)}\,\right),
$$
\begin{equation}
p_2=(8(\vert R^\ast_1 (r,\psi,\delta_0)\vert^2+
\vert R_0(r,\psi,\delta_0)\vert^2))/(1+\delta_0)^2=
\label{eq:bi}\end{equation}
$$
8(r^2-r\delta_0t+(\delta_0)^2/4+
\vert D_0(r,\psi,\delta_0)\vert)/(1+\delta_0)^2=
$$
$$
8(u-\delta_0\rho/2+(\delta_0)^2/4)/(1+\delta_0)^2+
$$
$$8(1+\delta_0)^{-2}\sqrt{u^2+u(\rho+1-(\delta_0)^2/2)+
(\delta_0/2)^2(\rho^2+\rho+(\delta_0/2)^2)}=
$$
$$
8(1+\delta_0)^{-2}u(2+(\rho+1-\delta_0\rho)/(2u)+O(1/u^2)),
$$
\begin{equation}\label{eq:cj}
q_1(r,\psi,\delta_0)=((1-\delta_0)/(1+\delta_0))^2,\,
q_2(r,\psi,\delta\sb 0)=
\end{equation}
$$
(4r/(1+\delta_0))^2=(16u)/(1+\delta_0)^2.
$$
In view of (91) in [\ref{r:cc}], (\ref{eq:be1}) and (\ref{eq:be3}),
\begin{equation}\label{eq:ca}
s=s_0(r,\psi)=\vert r\exp(i\psi)+1\vert\,/2=\sqrt{(r^2+1+2rcos(\psi))/4}=
\end{equation}
$$
\sqrt{(u+1+\rho)/4}\in(\max(\vert r-1\vert/2,\,\delta_0/4),\,(r+1)/2]
$$
and
$$
t=cos(\psi)=(4s^2-r^2-1)/(2r).
$$
In view of (3.1.68) in [\ref{r:bh}], (3.1.70) -- (3.1.71) in [\ref{r:bh}]
 and (\ref{eq:bg}),
$$
\vert R^\ast_{-1}(r,\psi,\delta_0)\vert^2=
r^2+(2+\delta_0)^2/4+r(2+\delta_0)\cos(\psi)=
$$
$$
u+(2+\delta_0)^2/4+\rho(2+\delta_0)/2,
$$
\begin{equation}
p_0=8(\vert R^\ast_{-1}(r,\psi,\delta_0)\vert^2+
\vert R_0(r,\psi,\delta_0)\vert^2)/(1+\delta_0)^2=
\label{eq:cd}\end{equation}
$$
8(u+(2+\delta_0)^2/4+\rho(2+\delta_0)/2)/(1+\delta_0)^2+$$
$$
8(1+\delta_0)^{-2}\sqrt{u^2+u(\rho+1-(\delta_0)^2/2)+
(\delta_0/2)^2(\rho^2+\rho+(\delta_0/2)^2)},
$$
\begin{equation}
q_0(r,\psi,\delta_0)(1+\delta_0)^2/16=(r^2+1+2rcos(\psi))=(u+1+\rho).
\label{eq:ce}\end{equation}
According to  Lemma 4.4 in [\ref{r:cc}], (\ref{eq:aa6}) and (\ref{eq:be3}),
\begin{equation}
\vert \eta_1^ \wedge (r,\psi,\delta_0)+\epsilon\vert<
\vert\eta_0^\wedge (r,\psi,\delta_0)+\epsilon\vert,
\label{eq:ch}\end{equation}
if $\epsilon^2=\epsilon$ and $z\in D_3.$
Therefore, according to (\ref{eq:bh}), (\ref{eq:cj})
and (\ref{eq:ch}),
\begin{equation}
(-1)^k(\partial/ \partial u)\vert\eta_k^\wedge(r,\psi,\delta_0)
\vert>0,
\label{eq:cf}\end{equation}
where $\frac13<\rho/2<\sqrt{u}=r,\,k^2=k.$
According
to  a) and c) of the Lemma 4.6 in [\ref{r:cc}],
and in view of (\ref{eq:aa6}) and  (\ref{eq:ca}),
\begin{equation}
\vert \eta_1^\wedge(r,\psi,\delta_0)-1\vert<
\vert\eta_0^\wedge (r,\psi,\delta_0)-1\vert,
\label{eq:ci}\end{equation}
if $z\in D_3.$
In view of (\ref{eq:bf}),
\begin{equation}
\vert D_0(r,\psi,\delta_0 )\vert^2 =
\label{eq:db}\end{equation}
$$u^2+u(\rho+1-(\delta_0)^2/2)+
(\delta_0/2)^2(\rho^2+\rho+(\delta_0/2)^2)=$$
$$(u+(\rho+1)/2-(\delta_0)^2/4)^2+
(\delta_0/2)^2(\rho^2+\rho+(\delta_0/2)^2)-$$
$$(((\rho+1)/2)^2-(\rho+1)(\delta_0)^2/4+
(\delta_0/2)^4)=$$
$$(u+(\rho+1)/2-(\delta_0)^2/4)^2+
(\delta_0/2)^2(\rho^2+2\rho+1)-(\rho+1)^2/4=$$
$$(u+(\rho+1)/2-(\delta_0)^2/4)^2-
(\rho+1)^2(1-(\delta_0)^2)/4.$$
Consequently,
\begin{equation}
\vert D_0(r,\psi,\delta_0 )\vert=
u+\frac{\rho+1}2-\frac{(\delta_0)^2}4+O(1/u),
\label{eq:db0}\end{equation}
where $u\ge1/4.$
Since $u\ge1/4>(\delta_0)^2/4,$ it follows that
$$u+(\rho+1)/2-(\delta_0)^2/4>\sqrt{1-(\delta_0)^2}(\rho+1)/2.$$
If $\rho=1,u=1/4$ then in view of (\ref{eq:db}),
$$
\vert D_0(r,\psi,\delta_0 )\vert^2 =
(5/4-(\delta_0)^2/4)^2-
(1-(\delta_0)^2)=$$
$$\left(\tau-5/4)^2\right)^2+4\tau-1,
$$
where $0<\tau=\frac{(\delta_0)^2}4<\frac1{100};$
moreover, in this case
$$(\partial/\partial\tau)\vert D_0(r,\psi,\delta_0 )\vert^2=
2\tau-5/2+4>0;$$ therefore if $\delta_0\le1/5,$ then
$$
\vert D_0(r,\psi,\delta_0 )\vert^2\bigg\vert_{u=1/4,\rho=1}\le
(1,24)^2-0,96=0,5776
$$
and
$$
\vert D_0(r,\psi,\delta_0 )\vert^2\bigg\vert_{u=1/4,\rho=1}\le0,76.
$$
In view of (\ref{eq:db}),
$$
1<(\partial/\partial u)\vert D_0(r,\psi,\delta_0 )\vert=
$$
$$
\sqrt{
\frac{(u+(\rho+1)/2-(\delta_0)^2/4)^2}
{(u+(\rho+1)/2-\frac{(\delta_0)^2}4)^2-
(\rho+1)^2(1-(\delta_0)^2)/4}}=1+O(1/u^2),$$
in view of (\ref{eq:bh}), (\ref{eq:bi}) and (\ref{eq:cd}),
\begin{equation}\label{eq:db2}
(\partial/\partial u)p_\epsilon=
8(2+O(1/u^2))/(1+\delta_0)^2,
\end{equation}
where $\epsilon^3=\epsilon,$
and $(\partial/\partial u)\vert D_0(r,\psi,\delta_0 )\vert$
decreases with increasing $u;$ consequently,
$$(\partial/{\partial u})^2
\vert D_0(r,\psi,\delta_0 )\vert<0,$$
if $u\ge1/4.$
In view of (\ref{eq:bh}), (\ref{eq:bi}) and (\ref{eq:cd}),
\begin{equation}
(\partial/\partial u)^2p_\epsilon=
(\partial/\partial u)^2
\vert D_0(r,\psi,\delta_0)\vert<0,
\label{eq:dc}\end{equation}
where $u\ge1/4,\,0<\delta_0<2/3<\rho,\,\epsilon^3=\epsilon.$
In view of (\ref{eq:bi}), (\ref{eq:cj}), (\ref{eq:db}) and (\ref{eq:db0}),
 if $\rho=1,\,u>1/4,\,0<\delta_0\le1/5,$ then
\begin{equation}\label{eq:da}
q_2((\partial/\partial u)p_2)/(\partial/\partial u)q_2-
p_2/2=\end{equation}
$$8u(1+(u+1-(\delta_0)^2/4))/
\vert D_0(r,\psi,\delta_0 )\vert)/(1+\delta_0)^2-$$
$$
4(u-\delta_0/2+(\delta_0)^2/4+
\vert D_0(r,\psi,\delta_0 )\vert)/(1+\delta_0)^2=
$$
$$
4(u+\delta_0/2-(\delta_0)^2/4)/(1+\delta_0)^2+
$$
$$
4((1+\delta_0)^2\vert D_0(r,\psi,\delta_0 )\vert)^{-1}
(2u^2+u(2-(\delta_0)^2/2)-$$
$$
4((1+\delta_0)^2\vert D_0(r,\psi,\delta_0 )\vert)^{-1}
(u^2+u(2-(\delta_0)^2/2)+
(\delta_0/2)^2(2+(\delta_0/2)^2))=
$$
$$
4(u+\delta_0/2-(\delta_0)^2/4)/(1+\delta_0)^2+
$$
$$
4(1+\delta_0)^2\vert D_0(r,\psi,\delta_0 )\vert)^{-1}
(u^2-(\delta_0/2)^2(2+(\delta_0/2)^2))>0,
$$
$$
q_2((\partial/\partial u)p_2)/(\partial/\partial u)q_2-
p_2=
\frac{8}
u(1+(u+1-(\delta_0)^2/4))/
\vert D_0(r,\psi,\delta_0)\vert)/(1+\delta_0)^2-
$$
$$
8(u-\delta_0/2+(\delta_0)^2/4+
\vert D_0(r,\psi,\delta_0 )\vert)/(1+\delta_0)^2=
$$
$$
8u(2+O(1/u^2))/(1+\delta_0)^2-
$$
$$
8(u-\frac{\delta_0}2+\frac{(\delta_0)^2}4+
u+1-\frac{(\delta_0)^2}4+O(1/u))/(1+\delta_0)^2=
$$
$$
-8(1-\frac{\delta_0}2+O(1/u))/(1+\delta_0)^2.
$$
In view of (\ref{eq:cd}), (\ref{eq:ce}), (\ref{eq:da}), (\ref{eq:db}),
(\ref{eq:db2}), (\ref{eq:db0}), if
 $\rho=1,\,u>1/4,\,0<\delta_0\le1/5,$ then
$$
(u+1)(\partial/\partial u)p_0-
p_0/2>
8(2u+2)/(1+\delta_0)^2-
$$
$$
4(u+(2+\delta_0)^2/4+(2+\delta_0)/2+u+1
-(\delta_0)^2/4)/(1+\delta_0)^2=$$
$$\frac8
(1/2+u-(3\delta_0)/4)/(1+\delta_0)^2>0,$$
\begin{equation}\label{eq:da3}
q_0((\partial/\partial u)p_0)/(\partial/\partial u)q_0-
p_0/2=
(u+2)(\partial/\partial u)p_0-
p_0/2>\end{equation}
$$(u+1)(\partial/\partial u)p_0-p_0/2>0,$$
\begin{equation}\label{eq:da4}
q_0(\partial/\partial u)p_0)/(\partial/\partial u)q_0-p_0=
8(u+2)(2+O(1/u^2))/(1+\delta_0)^2-
\end{equation}
$$
8(u+(2+\delta_0)^2/4+(2+\delta_0)/2+
u+1-(\delta_0)^2/4)/(1+\delta_0)^2=
$$
$$
8(4+O(1/u))/(1+\delta_0)^2-
(2+\delta_0)^2/4-(2+\delta_0)/2-1+(\delta_0)^2/4+
O(\frac1/u)=
$$
$$
8(1-(3/2)\delta_0+O(1/u))/(1+\delta_0)^2,
$$
where $u>1/4.$
In view of (\ref{eq:ce}), (\ref{eq:da3}) and (\ref{eq:dc}),
$$
(\partial/\partial u)
((
q_0(\partial/\partial u)p_0)/
(\partial/\partial u)q_0-p_0)(\partial/(\partial u)p_0+
(\partial/\partial u)q_0)=
$$
$$
(\partial/\partial u)(((u+2)(\partial/\partial u)p_0-
p_0)(\partial/\partial u)p_0)=
$$
$$
((\partial/\partial u)p_0)^2+
((u+2)(\partial/\partial u)^2p_0-
(\partial/\partial u)p_0)(\partial/\partial u)p_0+
$$
$$
(
(u+2)(\partial/\partial u)p_0-p_0)(\partial/\partial u)^2p_0=
$$
$$
((u+2)(\partial/\partial u)^2p_0)(\partial/\partial u)p_0+
((u+2)(\partial/(\partial u)p_0-p_0)
(\partial/\partial u)^2p_0=
$$
$$
(2(u+2)(\partial/\partial u)p_0-p_0)(\partial/\partial u)^2p_0<0.
$$
Therefore, according to (\ref{eq:da4}), (\ref{eq:db2}) and
 (\ref{eq:ce}),
\begin{equation}
\inf\{((u+2)(\partial/\partial u)p_0-p_0)(\partial/\partial u)p_0+
(\partial/\partial u)q_0\colon u\ge1/4\}=
\label{eq:gd}\end{equation}
$$
\lim\limits_{u\to+\infty}
((u(\partial/\partial u)p_0-p_0)(\partial/\partial u)p_0+
(\partial/\partial u)q_0)=
$$
$$
128(1-(3/2)\delta_0)/(1+\delta_0)^4+16/(1+\delta_0)^2)>0.
$$
According to the Lemma 4.17 in [\ref{r:cc}] and in view
 of (\ref{eq:da}), (\ref{eq:da3}), (\ref{eq:gd}),
\begin{equation}
(\partial/\partial u)\vert\eta_0(r,\psi,\delta_0)+\epsilon\vert^2>0,
\label{eq:ge}\end{equation}
where $\epsilon^2=1,\,u>1/4,$
\begin{equation}
(\partial/\partial u)\vert\eta_1(r,\psi,\delta_0)-1\vert^2<0,
\label{eq:gf}\end{equation}
where $u>1/4.$
The following Lemma describes the behavior of
 the value $h^\sim(\eta_k(r,\psi,\delta_0))$ with $k^2=k$ and $h^sim$
in (\ref{eq:ba}).

\bf Lemma 3. \it If $\Delta\ge5,$ then
\begin{equation}
(\partial/\partial u)(\vert h^\sim(\eta_0(r,\psi,\delta_0))\vert)>0,
\label{eq:gf1}\end{equation}
$$
(\partial/\partial u)(\vert h^\sim(\eta_1(r,\psi,\delta_0))\vert)<0,
$$
where $u\in(1/4,+\infty).$

\bf Proof. \rm The inequality (\ref{eq:gf1}) directly follows
from (\ref{eq:ch}), (\ref{eq:ge}) and (\ref{eq:ba}).
So, we must prove the inequality (\ref{eq:ba})
Clearly, if $\beta<1,\,u>1/4$ then
$$
(\partial/\partial u)(u^{3/4}+(3/4)\beta u^{-1/4})>0,
$$
We take
$$\beta=
(4/3)(\delta_0/2)^2
(2+(\delta_0)^2)/4)/(2-(\delta_0)^2/2).
$$
Then, clearly, $\beta<(\delta_0)^2=1/(\Delta)^2<1.$
Therefore,
 in view of (\ref{eq:bh}) and (\ref{eq:db}), if $\rho=1,$ then
$$
p_1u^{-1/4}=8(1+\delta_0)^{-2}\times
$$
$$
(u^{3/4}+(3/4)u^{-1/4}+
\sqrt{u^{3/2}+u^{1/4}(\rho+1-(\delta_0)^2/2)
(u^{3/4}+(3/4)\beta u^{-1/4})})
$$
increases together with increasing $u\in(1/4,+\infty),$
and, in view of (\ref{eq:cj}),
\begin{equation}\label{eq:gi}
\vert\eta_0(r,\psi,\delta_0)\vert^2u^{-1/4}\vert=
p_1u^{-1/4}/2+
\sqrt{(p_1u^{-1/4}/2)^2-q_1u^{-1/2}}
\end{equation}
increases together with increasing $u\in(1/4,+\infty).$

In view of (\ref{eq:cf}), (\ref{eq:cj}), (\ref{eq:gi}), (\ref{eq:ge})
and (\ref{eq:gf}), if $\Delta\ge5,$ then
$$
\vert\eta_1(r,\psi,\delta_0)\vert^{2(\Delta-1)}
\vert(\eta_1(r,\psi,\delta_0))^2-1\vert^2=
$$
$$
\vert\eta_1(r,\psi,\delta_0)\vert^{2(\Delta-5)}
\frac{(q_1)^4}{(\vert\eta_0(r,\psi,\delta_0)\vert^2u^{-1/4})^4}\times$$
$$\frac{16}{(1+\delta_0)^2}\vert\eta_0(r,\psi,\delta_0)+1\vert^{-2}
\vert\eta_1(r,\psi,\delta_0)-1\vert^2$$
decreases together with increasing $u\in(1/4,+\infty).\,\blacksquare$

Let $D$ is bounded domain in ${\mathbb C}$ or ${\mathfrak F}.$
and $D^\ast$ is closure of $D.$ Let
\begin{equation}
a^\sim_0(z)\,,\ldots\,,a^\sim_n(z)
\label{eq:dc1}\end{equation}
are the functions continuous on $D^\ast$ and analytic in $D.$
Let $a^\sim_n(z)=1$ for any $z\in D^\ast.$ Let
\begin{equation}
T(z,\lambda)=\sum\limits_{i=0}^na^\sim_i(z)\lambda^k.
\label{eq:dc2}\end{equation}
Let $s\in{\mathbb N},\,n_i\in{\mathbb N}-1,$ where $i=1,\,\ldots,\,s$
and
$\sum\limits_{i=1}^sn_i=n.$
We say that polynomial $T(z,\lambda)$ has
$(n_1,\,\ldots,\,n_s)$-disjoint system of roots on $D^\ast,$
 if for any $z\in D^\ast$ the set of all the roots $\lambda$
of the polynomial $T(z,\lambda)$ splits in $s$ klasses
${\mathfrak K}_1(z),\,\ldots,\,{\mathfrak K}_s(z)$
with  following properties:

 a) the sum of the multiplicities of the roots of the
klass ${\mathfrak K}_i$ is equal to $n_i$ for $i=1,\,\ldots,\,s;$

 b) if $i\in[1,s]\cap{\mathbb N},\,j\in(i,s]\cap{\mathbb N}$
and $n_in_j\ne0,$ then the absolute value of each roots
of the klass ${\mathfrak K}_i(z)$ is greater than absolute value
of the each roots of the klass ${\mathfrak K}_j(z).$

If the polynomial (\ref{eq:dc2}) has
$(n_1,\,\ldots,\,n_s)$-disjoint system of roots on $D^\ast,$
then for each $i=1,\,\ldots,\,s$ we denote by $\rho^\ast_{i,0}(z)$
and $\rho^\ast_{i,1}(z)$ respectively the maximal and minimal
absolute value of the roots of the klass ${\mathfrak K}_i(z).$

Let $D$ is bounded domain in ${\mathfrak F}$ such that
 $D^\ast\in D_3.$
Let
\begin{equation}
F^\wedge(z,\eta)=
\prod\limits_{i=1}^2(\theta_0(z)-h(\eta_{i-1}(r,\psi,\delta_0))),
\label{eq:bj3}\end{equation}
$$n=s=2,\,n_1=n_2=1,\,{\mathfrak K}_i(z)=\{h(\eta_{i-1}(r,\psi,\delta_0))\},$$
$$\rho_{i,0}=\rho_{i,1}=\vert h(\eta_{i-1}(r,\psi,\delta_0))\vert,$$
where $i=1,2.$

\bf Lemma 4. \it  The polynomial $F^\wedge(z,\eta)$ in (\ref{eq:bj3})
 has $(1,\,1)$-disjoint system of roots on $D^\ast.$

\bf Proof. \rm The assertion of the Lemma follows from (\ref{eq:ch})
 and (\ref{eq:ci}). $\blacksquare$

\bf Corollary. \it The map $(\ref{eq:bb})$ is injective
 for every  $z\in D^\ast;$ all the conditions of the Lemma 2 are fulfilled
 for the functions $f_0^\ast(z,\nu)$ from (\ref{eq:aa1}),
$\alpha^\ast(z,\nu)$ from (\ref{eq:ab}) and $\phi^\ast(z,\nu)$ from  (\ref{eq:ac})
 in every $z\in D^\ast;$ therefore for every $z\in D^\ast$ these functions
 are solutions of the difference equation of Poincar\'e type (\ref{eq:bc}),
 and the polynomial (\ref{eq:bd1})  coincides with
 characteristical polynomial of this equation. $\blacksquare$

Let for each $\nu\in{\mathbb N}-1$ are given continuous on $D^\ast$ functions
\begin{equation}\label{eq:dd0}
a_0(z;\nu),\,\ldots,\,a_n(z,\nu),
\end{equation}
which are analytic in $D.$

 Let $a_n(z:\nu)=1$
 for any $z\in D^\ast$ and any $\nu\in{\mathbb N}-1.$
We suppose that for any $i=1,\,\ldots,\,n-1$ the sequence
of functions  $a_i(z;\nu)$ converges to $a^\sim_i(z)$
 uniformly on $D^\ast,$  when $\nu\to\infty.$
Let us consider now the difference equation
\begin{equation}\label{eq:dd}
a_0(z;\nu)y(\nu+0)\,+\ldots\,+a_n(z;\nu)y(\nu+n)=0,
\end{equation}
i.e. we consider a difference equation of the Poincar\'e type,
 coefficients (\ref{eq:dd0}) of this equation are continuous on $D^\ast$
 and analytic in $D,$
and they uniformly converge to limit functions (\ref{eq:dc1}),
 when $\nu\to\infty.$

\bf Lemma 5. \it Let polynomial (\ref{eq:dc2}) has
$(n_1,\,\ldots,\,n_s)$-disjoint system of roots on $D^\ast.$
 Let $y(z,\,\nu)$ is a solution of the equation (\ref{eq:dd}),
and this solution is continuous on $D^\ast$ and analytic on $D.$
Let further $i\in[1,s]\cap{\mathbb Z}.$
Let us consider the set of all the $z\in D,$ for which the following
inequality holds
\begin{equation}
\limsup\limits_{\nu\in{\mathbb N},\,\nu\to\infty}
\vert y(z,\,\nu)\vert^{1/\nu})
<\rho_{i,1}(z);
\label{eq:de}\end{equation}
if this set has a limit point in $D,$ then the inequality (\ref{eq:de})
holds in $D^\ast.$

\bf Proof. \rm The proof may be found in [\ref{r:aa}]
 (Theorem 1 and its Corollary). $\blacksquare$

\bf Lemma 6. \it Let $D$ is bounded domain in ${\mathfrak F}$ such that
 $D^\ast\in D_3.$ Then
\begin{equation}
\limsup\limits_{\nu\in{\mathbb N},\,\nu\to\infty}
\vert f^\ast_0(z,\,\nu)\vert^{1/\nu})
<\rho_{1,1}(z)=\vert h^\sim(\eta_{0}(r,\psi,\delta_0))\vert
\label{eq:ha}\end{equation}
for any $z\in D^\ast.$

\bf Proof. \rm In view of (\ref{eq:aa6}), expanding the domain $D$,
 if necessary, we
can suppose that $\{(r,\,\phi)\colon\,r\in[2,\,3],\,\phi=0\}\in D.$
 Making use the same arguments, as in [\ref{r:bi}], Lemma 4.2.1,
 we see that the inequality (\ref{eq:ha})
 holds for any point $z=(r,\phi)\in\{r\in[2,\,3],\,\phi=0\}.$
 According to the Lemma 5, the inequality (\ref{eq:ha})
 holds for any $z\in D^\ast.\,\blacksquare$

For each prime $p\in{\mathbb N}$ let $v_p$ denotes the $p$-adic valuation
on $\Bbb Q.$

\bf Lemma 7.\it Let  $p\in{\mathbb N}+2 $ is a prime number,
$$
d\in{\mathbb N}-1,\,r\in{\mathbb N}-1,\,r<p.
$$
Then
$$
v_p((dp + r)!/((-p)^{d}d!\,r!)-1)\ge1.
$$

\bf Lemma 8. \it Let $p\in{\mathbb N}+2$ is a prime number,
$d\in{\mathbb N}-1,\,d_1\in{\mathbb N}-1$,
\begin {equation}\label{eq:zc}
r\in[0,p-1]\cap{\mathbb N},\,r_1\in[0,p-1]\cap{\mathbb N},\,d_1p+r_1\le dp+r.
\end{equation}
Then
\begin {equation}\label{eq:zd} v_p\left(
\binom{dp+r}{d_1p+r_1}\right)=v_p\left(\binom d{d_1}\right),
\end{equation}
if $r_1\le r,$
\begin {equation}\label{eq:ze} v_p\left(
\binom{dp+r}{d_1p+r_1}\left(\binom d{d_1}\binom r{r_1}\right)^{-1}-1\right)
\ge1,
\end{equation}
if $r_1\le r,$
\begin {equation}\label{eq:zf}
v_p\left(
\binom{dp+r}{d_1p+r_1}\right)=1+v_p\left((d-d_1)\binom d{d_1}\right),
\end{equation}
if $r<r_1,$
\begin {equation}\label{eq:zg}
v_p\left((-1)^{r_1-r-1}\binom{dp+r}{d_1p+r_1}
\binom{r_1}{r}(r_1-r)\left(p\binom{d}{d_1}(d-d_1)\right)^{-1}-1\right)\ge1,
\end{equation}
\bf Proof. \rm Clearly, $d_1\le d.$ If $r_1\le r,$ then
let $r_2=r-r_1,\,d_2=d-d_1.$ On the other hand,
if $r_1> r,$ then, in view of (\ref{eq:zc}), $d\ge d_1+1;$
therefore in this case we let
\begin{equation}\label{eq:zh}
r_2=p+r-r_1,\,d_2=d-d-1.
\end{equation}
Then $d=d_1+d_2,\,r=r_1+r_2,$
$$
\binom{dp+r}{d_1p+r_1}=(dp+r)!((d_1p+r_1)!(d_2p+r_2)!)^{-1}.
$$
Accordindg to the Lemma 7,
\begin {equation}\label{eq:zi}
v_p\left(\binom{dp+r}{d_1p+r_1}(-p)^{-d+d_1+d_2}
d_1!\,r_1!\,d_2!\,r_2!/(d!\,r!)-1\right)\ge1,
\end{equation}
\begin {equation}\label{eq:zaj}
v_p\left(\binom{dp+r}{d_1p+r_1}\right)=d-d_1-d_2+
\end{equation}
$$v_p(d!\,r!/(d_1!\,r_1!\,d_2!\,r_2!)).
$$
The equality (\ref{eq:zd}) and the inequality (\ref{eq:e})
diectly follow from (\ref{eq:zi}) and (\ref{eq:zaj}). If
\newpage
 the inequality $r<r_1$ holds, then in view of (\ref{eq:zh}) -- (\ref{eq:zaj}),
$$r_2!\prod\limits_{j=1}{r_1-r-1}(p+r-r_1+j)=(p-1)!,\,
v_p(r_2!(r_1-r-1)!(-1)^{r_1-r}-1)\ge1,$$
and (\ref{eq:zg}) holds.

\bf Corollary 1. \it Let  $p\in{\mathbb N} $ is a prime number,
$$
d\in{\mathbb N}-1,r\in{\mathbb N}-1,d_1\in{\mathbb N}-1,d_2\in{\mathbb N}-1,
r_1\in{\mathbb N}-1,r_2\in{\mathbb N}-1,
$$
$$
max(r_1,r_2)<p.
$$
Then
$$
p^{-d}(dp + r)!\in (-1)^d d!r!+ p\Bbb Z,
$$
$$
\binom{d_1+d_2)p+r_1+r_2}{d_1p+r_1}\in
\binom{d_1+d_2}{d_1}\binom{r_1+r_2}{r_1}+p{\mathbb Z}.
$$

\bf Proof. \rm This is direct corollary of the Lemma 7 and Lemma 8.
 See also Lemma 9 in [\ref{r:eh}].
$\blacksquare$

\bf Corolary 2. \it Let $p\in{\mathbb N}+2$ is a prime number,
$$
d\in{\mathbb N},\,r_1\in{\mathbb N},\,r_1<p,\,d_1\in{\mathbb N}-1,\,d_1<d.
$$
Then
\begin {equation}\label{eq:zf1}
v_p\left(\binom{dp}{d_1p + r_1}\left(
 d\binom{d-1}{d_1}{\binom p{r_1}}\right)^{-1}+1\right)\ge1
\end{equation}

\bf Proof. \rm
Since,
$$d\binom{d-1}{d_1}=(d-d_1)\binom d{d_1},
\,v_p\left(\binom p{r_1}r_1/p-(-1)^{r_1}\right)\ge1,$$
the equality (\ref{eq:zf1}) directly follows from (\ref{eq:zg}). $\blacksquare$

\bf Corolary 3. \it Let $p\in{\mathbb N}+2$ is a prime number,
$$
d\in{\mathbb N},\,r_1\in{\mathbb N},
\,r_1<p,\,d^\sim\in{\mathbb N}-1,\,d^\sim<d.
$$
Then
$$
\binom{dp}{d_1p + r_1}\in
 d\binom{d-1}{d^\sim}\binom p{r_1}+p^2{\mathbb Z}.
$$
\bf Proof. \rm  This is a corollary
of the Corrolary 2. See also Lemma 10 in [\ref{r:eh}] .$\blacksquare$

Let let $p$ be prime in $(2,\,+\infty),$
let $K$ be a finite extension of ${\mathbb Q}$  
let ${\mathfrak p}$ be a prime ideal in ${\mathbb Z}_K$
 and $p\in{\mathfrak p},$ let $f$ be the degree of ${\mathfrak p},$
let $(p)={\mathfrak p}^e{\mathfrak b},$ with entire ideal  ${\mathfrak b}$ 
 not contained in ${\mathfrak p},$ let $v_{\mathfrak p}$
be additive ${\mathfrak p}$-valuation, which prolongs $v_p;$ so,
 if $\pi$ is a ${\mathfrak p}$-prime number, then 
$v_{\mathfrak p}(\pi)=1/e.$ If $f$ is the degree of the ideal ${\mathfrak p}$
then 
\begin{equation}\label{eq:z7a}
v_{\mathfrak p}\left(w^{p^\beta}-w\right)\ge1,
\end{equation}
where $\beta\in{\mathbb N}f,\,w\in K$ and
$$
v_{\mathfrak p}(w)\ge0.
$$
In viw of (\ref{eq:z7a}), (\ref{eq:ab}),  and (\ref{eq:i}),
$$
v_{\mathfrak p}(\alpha^\ast(z;\,p^\beta l)-\alpha^\ast(z;\,l))>1/e,
$$
if $\beta\in{\mathbb N}f,\,\theta_0(z)\in K$ and
 $v_{\mathfrak p}(\theta_0(z))\ge0.$
In view of (\ref{eq:ac}),
\begin{equation}\label{eq:r7d}
\phi^\ast(z;\nu)=(-\theta_0(z))^\nu
\sum\limits_{k=0}^{\nu\Delta}\alpha^\ast_{\nu,k}(\theta_0(z))^{k}
\sum \limits_{\tau=1}^{\nu+k}
((\theta_0(z))^{-\tau}/\tau))=
\end{equation}
$$
(-\theta_0(z))^\nu\sum \limits_{\tau=1}^\nu((\theta_0(z))^{-\tau}
\alpha^\ast(z;\nu)/\tau+
$$
$$
(-\theta_0(z))^\nu
\sum\limits_{k=0}^{\nu\Delta}\alpha^\ast_{\nu,k}(\theta_0(z))^{k}
\sum \limits_{\tau=1+\nu}^{\nu+k}
((\theta_0(z))^{-\tau}/\tau))=
$$
$$(-1)^\nu\sum\limits_{\tau=1}^{\nu(\Delta+1)}\frac1\tau\sum
\limits_{k=\max(0,\,\tau-\nu)}^{\nu\Delta}
\alpha^\ast_{\nu,k}(\theta_0(z))^{\nu-\tau+k};$$
therefore, if $\nu=p^\beta l,\,f=1,\,\beta\in{\mathbb N}f,\,p>l(\Delta+1)
,\,\theta_0(z)\in K$ and $v_{\mathfrak p}(\theta_0(z))\ge0,$
 then, according to the Lemma 2,
\begin{equation}\label{eq:z7e}1-\beta\le\end{equation}
$$v_{\mathfrak p}\left(\phi^\ast(z;\nu)-
\sum\limits\Sb\eta\in[1,\,\Delta+1]\cap {\mathbb Z}\\
k\in[p^\beta(\eta-l),\,p^\beta l\Delta]\cap {\mathbb Z}\\
k\ge0,\, v_{\mathfrak p}(k)>0\endSb
\frac{(-1)^{pl}}{p^\beta\eta}
(\theta_0(z))^{p^\beta(l-\eta)+k}\alpha^\ast_{\nu,k}\right),$$
\begin{equation}\label{eq:z7f}1/e-\beta\le\end{equation}
$$v_{\mathfrak p}\left(\phi^\ast(z;\nu)-
\sum\limits\Sb\eta\in[1,\,\Delta+1]\cap {\mathbb Z}\\
k\in[p^{\beta-1}(\eta-l),\,p^{\beta-1}l\Delta]\cap {\mathbb Z}\\
k\ge0\endSb
\frac{(-1)^{pl}p}{p^\beta\eta}
(\theta_0(z))^{p^{\beta-1}(l-\eta)+k}\alpha^\ast_{\nu/p,k}\right).
$$
We make the pass (\ref{eq:z7e}) $\to$ (\ref{eq:z7f}) $\beta$ times
and obtain the inequality
\begin{equation}\label{eq:z7g}1/e-\beta\le\end{equation}
$$v_{\mathfrak p}\left(\phi^\ast(z;\,p^\beta l)-p^{-\beta}\phi^\ast(z;\, l)
\right),
$$
where
$$\{l,\,\beta\}\subset{\mathbb N},\,p>l(\Delta+1),\,p\in\mathfrak p$$
and $\mathfrak p$ is ideal of the first degree.

\bf Lemma 9. \it If $m\in{\mathbb N}+1,\,K={\mathbb Q}[\exp(2pi/m)],$
$$
\alpha^\ast(z;\,l_1)\phi^\ast(z;\,l_2))\ne0
$$
for some $z\in K\diagdown\{0\},\,l_1\in{\mathbb N},\,l_2\in{\mathbb N},$
then for any $l\in{\mathbb N}$ the sequenses
\begin{equation}\label{eq:z7i}
\alpha^\ast(z;\,\nu),\,\phi^\ast(z;\,\nu),
\end{equation}
where $\nu\in l+{\mathbb N}$
form a linear independent system over $K.$

\bf Proof. \rm There exists $d^\ast\in{\mathbb N}$ such that
$$d^\ast z\in{\mathbb Z}_K,\,d^\ast z\alpha^\ast(z;\,l_1)\in{\mathbb Z}_K,\,
d^\ast z\phi^\ast(z;\,l_2)\in{\mathbb Z}_K.$$
Let a prime $p\in{\mathbb N}m+1$ satisfies to the inequality
$$
p>\vert \Nm_{K/{\mathbb Q}}(d^\ast z\alpha^\ast(z;\,\nu))\vert+
\vert \Nm_{K/{\mathbb Q}}(d^\ast z\phi^\ast(z;\,\nu))+
$$
$$\vert \Nm_{K/{\mathbb Q}}(d^\ast z)\vert+
\vert \Nm_{K/{\mathbb Q}}(d^\ast)+(\Delta+1)(l_1+l_2).$$
Let $\mathfrak p$ is a prime ideal containing $p.$
Then
$$
v_{\mathfrak p}\left(\alpha^\ast(z;\,l_1)\right)=
v_{\mathfrak p}\left(\phi^\ast(z;\,l_2)\right)=0,
$$
and, in view of (\ref{eq:z7g}), 
$$
v_{\mathfrak p}\left(\phi^\ast(z;\,p^\beta l_2)\right)=-\beta,
$$
but
$$
v_{\mathfrak p}\left(\alpha^\ast(z;\,p^\beta l_1)\right)=0.
$$
with $\beta\in{\mathbb N}\,\blacksquare.$

Let $m\in{\mathbb N},\,k\in{\mathbb Z},\,2\le2\vert k\vert<m,$
and let $m$ and $k$ have no common divisor
with exeption $\pm1.$  Let further
$K_m={{\mathbb Q}[\exp(2\pi i/m)]}$
is a cyclotomic field, ${\mathbb Z}_{K_m}$ is the ring of all the integers
of the field $K_m.$ 

\bf Lemma 10. \it Let $\Delta\in\{5,\,7\}.$
 In correspondece with (\ref{eq:aa3}),
(\ref{eq:aa4}) and (\ref{eq:aa6}),
let $z=\left(1/(2\cos(k\pi i/m),k\pi i/m-\pi\right),$
where $\vert k\vert<m/2,\, (\vert k\vert,m)=1.$ 

Then for each $l\in{\mathbb N}$ the two sequences (\ref{eq:z7i})
form a linear independent system over $\mathbb C.$

\bf Proof. \rm We check the fulfilment of the conditions of theLemma 9.

 Let ${\mathfrak M}={\mathbb N}\diagdown\{1,\,2,\,6\}$ and  
${\mathfrak M}_0=\{m\in{\mathfrak M}\colon \Lambda_0(m)=0\}.$  
According to the condition of
the Lemma, $\theta_0(z)=-1/(1+\exp(2i\pi/m)$ with $m\in{\mathfrak M}.$
If $m\in{\mathfrak M}$ and $\phi(m)>\Delta$, then, in view of
(\ref{eq:ab}) and (\ref{eq:i}), $\alpha^\ast(z;1)\ne0,$ because the numbers
$(1+\exp(2i\pi/m)^k,$ where $k=0,\,\ldots,\,\phi(m)-1,$ form a basis of
the field $K_m.$ Let $\Delta=p\in2{\mathbb N}+1,$
 where $p$ is a prime, 
 $\mathfrak p$ is a prime ideal containing $p,$ 
and, as before,
 let $(p)={\mathfrak b}{\mathfrak p}^e,\,1_{K_m}\in{\mathfrak b}\,+\,
{\mathfrak p}.$
Then
\begin{equation}\label{eq:z13}
\binom{2p-1}p\binom p{p-1}\equiv p \mod p^2,\,
v_{\mathfrak p}\left(\binom{p+k}{1+k}\binom pk\right)=2,
\end{equation}
where $k=1,\,\ldots,\,p-2,$
\begin{equation}\label{eq:z14}
\binom p1\binom p0=p,\,\binom{2p}{p+1}\binom pp\equiv 2p \mod p^2.
\end{equation}
If $m\in{\mathfrak M}$ and $(m,p)=1,$ or, if $m\in{\mathfrak M}_0,$
then, according to the Lemma 1,
\begin{equation}\label{eq:z16}
(1+\exp(2i\pi/m),p)=(1)
\end{equation}
and, according to the Lemmata 7 and 8,
\begin{equation}\label{eq:z17}
\alpha^\ast(z;1)/(p\theta_0(z))\equiv 1+(\theta_0(z))^{p-1}-
2(\theta_0(z))^p\equiv
\end{equation}
$$ 1+(\exp(2i\pi/m)+3)/(1+\exp(2pi\pi/m))\equiv$$
$$(\exp(2ip\pi/m)+\exp(2i\pi/m)+4)/(1+\exp(2pi\pi/m)) \mod p. $$
If $m=q^\alpha$ with $\alpha\in{\mathbb N}$ and prime $q$
and there exists $l$ in $\{0,\,\ldots,\,\phi(m)-1\}$ such that
 $p\equiv l \mod (m),$ then
\begin{equation}\label{eq:y1}
\exp(2ip\pi/m)+\exp(2i\pi/m)+4\not\equiv0 \mod p.
\end{equation}
If $m=2q^\alpha$ with odd prime $q$ and $\alpha\in{\mathbb N,}$ and
 there exist $l$
in $\{0,\,\ldots,\,\phi(m/2)-1\}$ such that
 $p\equiv 2l\mod (m/2),$ then (\ref{eq:y1}) holds.

If $p=5,$ then $\{3,\,4,\,5\,\,8,\,10,\,12\}=
\{m\in\mathfrak M\colon\phi(m)\le p\}.$

 If $m=3,\,4,\,5,\,8,\,10$ then, clearly, (\ref{eq:y1}) holds.

If $m=12,$ then $1,\,\exp(i\pi/2),\,\exp(2i\pi/3),\,\exp(i\pi/6),\,$
form a entire basis of $K_{12},$ 
$\exp(5i\pi/6)=\exp(i\pi/2)-\exp(i\pi/6),\,$ and (\ref{eq:y1}) holds.

If $p=7$ then $\{3,\,4,\,5\,7,\,\,8,\,9,\,10,\,12,\,14,\,18\}=
\{m\in\mathfrak M\colon\phi(m)\le p.$

 If $m=3,4,\,\,5,\,7,\,9,\,14,$ then, clearly, (\ref{eq:y1}) holds.

If $m=8,$ then $\exp(7i\pi/4)=-\exp(3i\pi/4)$ and (\ref{eq:y1}) holds.

If $m=12,$ then $1,\,\exp(i\pi/2),\,\exp(2i\pi/3),\,\exp(i\pi/6),\,$
form a entire basis of $K_{12},$ 
$\exp(7i\pi/6)=-\exp(i\pi/6),\,$ and (\ref{eq:y1}) holds.

If $m=18,$ then
$$\exp(7i\pi/9)=-\exp(-2i\pi/9)=\exp(4i\pi/9)+\exp(10i\pi/9),\,$$
 and (\ref{eq:y1}) holds.

The coefficient at $(\theta_0(z))^0$
in the expression (\ref{eq:ac}) of $\phi^\ast(z;\nu)$ 
is equal to
$$\sum\limits_{k=0}^{\nu\Delta}(-1)^\nu\alpha^\ast_{\nu,k}/(\nu+k)$$
and, if $\Delta=p,\,\nu=1,$ then in view (\ref{eq:z13}) -- (\ref{eq:z14}),
 the value of $v_p$ on this coefficient is equal to $0.$
 Therefore,
 if $m\in{\mathfrak M}$ and $\phi(m)>p=\Delta$, then $\phi^\ast(z;1)\ne0.$  

If $m\in{\mathfrak M}\diagdown{\mathfrak M}_0,$ and $m\equiv0 \mod p$
 then $m=2p^\alpha,$ where $\alpha\in{\mathbb N}.$
 According to the Lemma 1, ${\mathfrak p}=(1+\exp(2i\pi/m)$ is a prime ideal
in $K_m,$ and, furthermore, ${\mathfrak p}^{\phi(m)}=(p).$
 Let $v_{\mathfrak p}.$ is the ${\mathfrak p}$-adic valuation,
 which prolongs the valuation $v_p.$
Clearly, $v_{\mathfrak p}(1+\exp(2i\pi/m)=1/\phi(m),\,
v_{\mathfrak p}(\theta_0(z))=-1/\phi(m)$
In view of (\ref{eq:ac}) with $\nu=1$,
for the summands of the sum
$$
\sum\limits_{k=1}^{\nu\Delta}\alpha^\ast_{\nu,k}(\theta_0(z))^{1+k}
\sum \limits_{\tau=2}^{1+k}
((\theta_0(z))^{-\tau}/\tau))
$$
we have the inequality
$$v_{\mathfrak p}((\theta_0(z))^{\Delta+k-y\alpha^\ast_{\nu,k}/y}\ge
-(k-1)/\phi(m)+2-v_{\mathfrak p}(\tau)\ge-(p-3)/\phi(m)+2,$$
if $k=1,\,\ldots,\,p-2,$ because in this case $\tau\in[2,\,p-1],$
$$v_p((\theta_0(z))^{\Delta+k-y\alpha^\ast_{\nu,k}/y}\ge-(k)/\phi(m)+
1-v_p(\tau)\ge-(p-1)/\phi(m),$$ 
where $k\in\{q-1,\,q\},$ and the equality reaches only for $k=\tau=p;$
on the other hand,
 $v_{\mathfrak p}(\alpha^\ast(z;1))\ge1-(p+1)/\phi(m)\ge
-2/(p-1)\ge-2/\phi(m).$  
So, if $p\ge5,$ then $v_{\mathfrak p}(\phi^\ast(z;1))=-(p-1)/\phi(m).$
If $m\in{\mathfrak M}\diagdown{\mathfrak M_0,}$ then $m=2q^\alpha,$
with prime $q,$ according to the Lemma 1, ${\mathfrak l}=(1+\exp(2i\pi/m))$
is a prime ideal in $K_m,$ and ${\mathfrak l}^{\phi(m)}=(q).$
Therefore in this case $v_{\mathfrak p}(\theta_0(z)=0$
If $m\in{\mathfrak M}_0,$ then, according to the Lemma 1,
 $v_{\mathfrak p}(\theta_0(z)=0.$
According to (\ref{eq:ac}), in both last cases,
$$v_{\mathfrak p}(\phi^\ast(z;1))+
\alpha_{\nu,p-1}/p+\theta_0(z)\alpha_{\nu,p}/p)\ge1.$$
In view of (\ref{eq:z13}),\,(\ref{eq:z14}), 
$$v_{\mathfrak p}(\alpha_{\nu,p-1}/p+\theta_0(z)\alpha_{\nu,p}/p)=$$
$$v_{\mathfrak p}(exp(2i\pi/m)-1)/(exp(2i\pi/m)+1)).$$
If $p=5$ and $m\in\{3,\,4,\,5,\,7,\,8,\,9,\,10\}$
then, clearly,
\begin{equation}\label{eq:y5}
v_{\mathfrak p}((exp(2i\pi/m)-1))\le1/4.
\end{equation}
If $p=5$ and $m=12,$ then $\Nm_{K_{12}}((exp(i\pi/6)-1))=3$
and (\ref{eq:y5}) holds. 

If $p=7,$ and $m\in\{3,\,4,\,5,\,7,\,8,\,9,\,10,\,12,\,14,\,18\},$
 then $$v_{\mathfrak p}((exp(2i\pi/m)-1))\le1/6.$$
 $\blacksquare$

\bf Lemma 11. \it Let are fulfilled all the conditions of the Lemma 10.
Then
\begin{equation}
\limsup\limits_{\nu\in{\mathbb N},\,\nu\to\infty}
\left(\vert f^\ast_0(z,\,\nu)\vert^{1/\nu}
=\rho_{2,1}(z)\bigg\vert_{\theta_0(z)=-1/(1+\exp(2ik\pi/m)}\right)=
\label{eq:hd}\end{equation}
$$\vert h^\sim(\eta_{1}
(1/(2\cos(k\pi i/m)),\,k\pi i/m,\delta_0))\vert,$$
where $h^\sim(\eta)$ is defined in (\ref{eq:ba}).

\bf Proof. \rm According to the Lemma 2, (\ref{eq:aa1}) and  Lemma 10,
$f^\ast_0(z,\,\nu)$ is a nonzero solution of the Poincar\'e type
difference equation (\ref{eq:bc}). According to the Perron's theorem
and Lemma 5, the equality (\ref{eq:hd}) holds. $\blacksquare$

Let $K/{\mathbb Q}$ be the finite extension of the field ${\mathbb Q},$
$$[K:{\mathbb Q}] = d.$$ Let the field $K$ has $r_1$ real places
 and $r_2$ complex places. Each such place is the monomorphism
 of the field $K$ in the field ${\mathbb R},$ if a place is real, or in the
field ${\mathbb C},$ if a place is not real; we will denote these monomorphisms
respectively by $\sigma_1\,,\ldots\sigma_{r_1+r_2}.$
Then $d=r_1+2r_2.$ Let ${\mathfrak B}$ be the fixed integer basis
$$
\omega_1\,,\ldots\,,\omega_d
$$
of the field $K$ over ${\mathbb Q}.$ Clearly,$K$ is an algebra over
 ${\mathbb Q}.$
With extension of the ground field from ${\mathbb Q}$
to ${\mathbb R}$  appears an isomorphism
of the algebra ${\mathfrak K} = K\otimes{\mathbb R}$ onto direct sum
$$
\underbrace{\mathbb R\oplus\ldots\oplus\mathbb R}_{\text {$r_1$ times}}
\oplus\underbrace{\mathbb C\oplus\ldots\oplus\mathbb C}_
{\text {$r\sb 2$ times}}
$$
of $r_1$ copies of the field ${\mathbb R}$ and
$r_2$ copies of the field ${\mathbb C}.$
We identify by means of this isomorphism the aIgebra ${\mathfrak K}$
with the specified direct sum. We denote below by
 $\pi_j,$ where $j = 1\,,\ldots\,,r_1+r_2,$ the projection
 of ${\mathfrak K}$ onto its $j-$th direct summand
 and also the extension of this projection
onto all kinds of matrices which have all the elements in ${\mathfrak K}.$
 So, $\pi_j({\mathfrak K})={\mathbb R}$
 for $j = 1\,,\ldots\,, r_1$ and $\pi_j({\mathfrak K})={\mathbb C}$ for
 $j = r_1+1\,,\ldots\,,r_1+r_2.$
 Further by ${\mathfrak i}_{\mathfrak K}$ we denote the embedding
 of ${\mathbb R}$ in ${\mathfrak K}$ in diagonal way and also
 the extension of this embedding
 onto all kinds of the real matrices. So, ${\mathbb R}$ is imbedded by means
 of ${\mathfrak i}_{\mathfrak K}$ in ${\mathfrak K}$ in diagonal way.
 Each element $Z\in{{\mathfrak K}}$ has a unique representation in the form:
$$
Z=\left(\matrix z_1\\\vdots\\ z_{r_1+r_2}\\\ \\\overline{z_{r_1+1}}
\\\vdots\\\overline{z_{r_1+r_2}}\endmatrix\right),
$$
with $z_j=\pi_j(Z)\in{\mathbb R}$ for any $j=1\,,\ldots\,,r_1$ and with
 $z_j=\pi_j(Z)\in\mathbb C$ for

\noindent any $j=r_1+1\,,\ldots\,,r_1+r_2.$
Further by $\Tr_{\mathfrak K} (Z)$ we denote
the sum
$$\sum\limits_{j=1}^{r_1}z_j+\sum\limits_{j=r_1+1}^{r_1+r_2}2\Re(z_j)=$$
$$\sum\limits_{j=1}^{r_1}\pi_j(Z)+\sum\limits_{j=r_1+1}^{r_1+r_2}
 2\Re(\pi_j (Z)),$$
and by $q_\infty^{({\mathfrak K})} (Z)$ we denote
the value
$$\max (\vert z_1\vert\,,\ldots\,,\vert z_{r_1+r_2}\vert)=$$
$$\max(\vert\pi_1(Z)\vert\,,\ldots\,,\vert\pi_{r_1+r_2}(Z)\vert).$$
Clearly,
$$q_\infty^{({\mathfrak K})}(Z_1Z_2)\le q_\infty^{({\mathfrak K})}(Z_1)
q_\infty^{({\mathfrak K})}(Z_2),$$
$$q_\infty^{({\mathfrak K})}(Z_1+Z_2)\le q_\infty^{({\mathfrak K})}(Z\sb1)
+q_\infty^{({\mathfrak K})}(Z_2),$$
$$q_\infty^{(\mathfrak K)}({\mathfrak i}_{\mathfrak K}(\lambda)Z)=
\vert\lambda\vert q_\infty^{({\mathfrak K})}(Z)$$
for any $Z_1\in{\mathfrak K},\,Z_2\in {\mathfrak K},\,Z\in{\mathfrak K}$
and $\lambda\in{\mathbb R}.$
The natural extension of the norm $q_\infty^{({\mathfrak K})}$
on the set of all the matrices, which have all the elements in $\mathfrak K$
(i.e. the  maximum of the norm $q_\infty^{({\mathfrak K})}$
of all the elements of the matrix) also will be denoted
by $q_\infty^{({\mathfrak K})}.$
If $$Z=\left(\matrix z_1\\\vdots\\z_{d}\endmatrix\right)\in K,$$
then
$$z\sb j=\sigma_j (Z),$$
where $j = 1\,,\ldots\,,r_1+r_2,$
$$z_{r_1+r_2+j}=\overline {\sigma_{r_1+j}(Z)},$$
where $j = 1\,,\ldots\,,r\sb 2.$
In particular,
$$
\omega_k=\left(\matrix\sigma_1(\omega_k)\\\vdots\\\
\sigma_{r_1+r_2 }(\omega_k)\\\ \\\overline{\sigma_{r\sb1+1}(\omega_k)}
\\\vdots\\\overline{\sigma_{r_1+r_2}(\omega_k)}
\endmatrix\right),
$$
As usually, the ring of all the integer elements of the field $K$
 will be denoted by ${\mathbb Z}_K.$
 The ring ${\mathbb Z}_K$ is embedded in the ring ${\mathfrak K}$ as discrete
lattice. Moreover, if $Z\in{\mathbb Z}_K\backslash \{0\},$ then
$$\left (\prod\limits_{i=1}^{r_1}\vert\sigma_j(Z)\vert\right)
\prod\limits_{i=1}^{r_2}\vert\sigma_{r_1+i}(Z)\vert^2 =
\vert\Nm_{K/{\mathbb Q}}(Z)\vert \in{\mathbb N}$$
and therefore
$
q_\infty^{({\mathfrak K})}(Z)\ge1.
$
for any $Z\in {\mathbb Z}_K\backslash \{0\}.$ The elements of ${\mathbb Z}_K$
we name below by $K$-integers.
For each $Z\in{\mathfrak K}$ let
$$\Vert{\mathbb Z}\Vert_K
=\inf\limits_{W\in{\mathbb Z}_K}\{q^{({\mathfrak K})}_\infty(Z-W)\}. $$
Let $\{m,\,n\}\subset \Bbb N,$
$$
a_{i,k}\in{\mathfrak K}
$$
 for $i = 1\,,\ldots\,,m,\ k = 1\,,\ldots\,,n,$
$$
\alpha_j^\wedge(\nu)\in{\mathbb Z}_K,
$$
where $j=1\,,\ldots \,,m+n$ and $\nu\in\Bbb N.$
Let there are $\gamma_0,r^\wedge_1\ge 1,\,\ldots,\,r^\wedge_m\ge 1$
such that
$$
q^{(\mathfrak K)}_\infty(\alpha_i(\nu))<
\gamma_0(r^\wedge_i)^\nu$$
where $i=1\,,\ldots \,,m$ and $\nu\in\Bbb N.$
Let
$$
y_k(\nu)=-\alpha^\wedge_{m+k}(\nu)+
\sum\limits_{i=1}^ma_{i,k}\alpha_i^\wedge(\nu)$$
where $k=1\,,\ldots \,,n$ and $\nu\in{\mathbb N}.$
If
$X=\left(\matrix Z_1\\\vdots\\Z_n\endmatrix\right)\in{\mathfrak K}^n,$
 then let
$$
y^\wedge (X)=y^\wedge(X,\nu)=\sum\limits_{k=1}^ny^\wedge_k(\nu)Z_k
$$
for $\nu\in{\mathbb N},$
let
$$
 \phi_i(X)=\sum\limits_{k=1}^n a_{i,k} Z_k
$$
for $i=1\,,\ldots \,,m,$
and let
$$
\alpha^\wedge_0(X,\nu)=\sum\limits_{k=1}^n\alpha^\wedge_{m+k}(\nu)Z_k
$$
for $\nu\in\mathbb N.$
Clearly,
$$
y^\wedge (X,\nu)=-\alpha_0^\wedge (X,\nu)+
\sum\limits_{i=1}^m\alpha_i^\wedge(\nu)\phi_i(X)
$$
for $X\in{\mathfrak K}^n$ and $\nu\in\Bbb N,$
$$
 \alpha^\wedge_0(X,\nu)\in{\mathbb Z_K}
$$
for $X \in ({\mathbb Z}_K)^n$ and $\nu\in\Bbb N.$

\bf Lemma 12. \it Let $\{l,\,n\}\subset{\mathbb N},\,\gamma_1>0,\,
\gamma_2>\frac 12,R_1\ge R_2>1,$
$$
\alpha_i=(\log(r^\wedge_iR_1/R_2))/\log(R_2) ,
$$
where $i=1\,,\ldots \,,m,$ let
$
X\in ({\mathbb Z}_K)^n\backslash\{(0)\},$
$$
\gamma_3=\gamma_1(R_1)^{(-\log(2\gamma_2R_2))/\log(R_2)} ,
\gamma_4=\gamma_3
\left(
\sum\limits_{i=1}^m\gamma_0(r_i^\wedge)^{(log(2\gamma_2))/\log(R_2)+l}
\right)^{-1}$$
and let for each $\nu\in{\mathbb N} - 1$ hold the inequalities
$$
\gamma_1(R_1)^{-\nu}q^{({\mathfrak K})}_\infty(X)\le
\sup\{q^{({\mathfrak K})}_\infty(y^\wedge(X,\kappa))\colon\kappa=
\nu,\,\ldots,\,\nu+l-1\},$$
$$
q\sp{({\mathfrak K})}_\infty(y^\wedge(X,\nu))\le\gamma_2(R_2)^{-\nu}
q^{(\mathfrak K)}_\infty(X)
$$
Then
$$
\sup\{\Vert\phi_i(X)\Vert_K(q^{({\mathfrak K})}_\infty (X))^{\alpha_i}
\colon i=1,\,\ldots, \,m\} \ge\gamma_4.
$$

\bf Proof. \rm Proof may be found in [\ref{r:cj}], Theorem 2.3.1. $\blacksquare$

\bf Corollary. \it Let
$
a\in {\mathfrak K},
$
\begin{equation}
\alpha_1^\wedge(\nu)\in{\mathbb Z}_K,\,\alpha_2^\wedge(\nu)\in{\mathbb Z}_K,
y(\nu)=-\alpha^\wedge_2(\nu)+a\alpha_1^\wedge (\nu)
\label{eq:ec}\end{equation}
where $\nu\in\mathbb N.$
Let there are $\gamma_0,r^\wedge_1\ge1$ such that
$$
q^{({\mathfrak K})}_\infty(\alpha_1(\nu))<\gamma_0(r^\wedge_1)^\nu,
$$
where $\nu\in{\mathbb N}.$
 Let $l\in{\mathbb N},\,\gamma_1>0,\,
\gamma_2>\frac 12,R_1\ge R_2>1,$
$$
\alpha_1=(\log(r_1^\wedge R_1/R_2))/\log(R_2),\,
\gamma_3=\gamma_1(R_1)^{(-\log(2\gamma_2R_2))/\log(R_2)},
$$
$$
\gamma_4=\gamma_3\left (
 \gamma_0(r^\wedge_1)^{(log(2\gamma_2))/\log(R_2)+l}\right)^{-1},
$$
$X\in{\mathbb Z}_K$ and let for each $\nu\in{\mathbb N}-1$
 hold the inequalities
$$
\gamma_1(R_1)^{-\nu}q^{({\mathfrak K})}_\infty(X)\le
\sup\{q^{({\mathfrak K})}_\infty(y_1(\kappa)X)\colon
\kappa=\nu\,,\ldots\,,\nu+l-1))\},
$$
$$
q^{({\mathfrak K})}_\infty(y(\nu)X)\le\gamma_2
(R_2)^{-\nu}q^{({\mathfrak K})}_\infty(X)
$$
Then
\begin{equation}\label{eq:fc}
\Vert aX\Vert_K(q^{({\mathfrak K})}_\infty (X))^{\alpha}
 \ge\gamma_4.
\end{equation}

\bf Proof. \rm This Corrolary is the Lemma 12 for $m=n=1.$ $\blacksquare$

Let $B\in{\mathbb N},\,
D^\ast(B)=\inf\{q\in{\mathbb N}\colon d/\kappa \in{\mathbb N},\,
\kappa \in{\mathbb N},\,\kappa\le B\}.$
It is known that $$D^\ast(B)=\exp(B+O(B/\log(B)).$$ Let
$d^\ast_0(\Delta,\nu)=D^\ast(\nu(\Delta+1)).$ Then
\begin{equation}
d^\ast_0(\Delta,\nu)=\exp(\nu(\Delta+1)+O(\nu/\log(\nu))),
\label{eq:fd}\end{equation}
when $\nu\to\infty.$

Probably G.V. Chudnovsky was the first man, who discovered, that
the numbers (\ref{eq:i}) have a  great common divisor; Hata ([\ref{r:eb}]) in
details studied this effect. Therefore I name the mentioned common divisor
 by Chudnovsky-Hata's multiplier and denote it
by $d^\ast_1(\Delta,\nu).$ According to the Hata's results,
\begin{equation}
\log(d^\ast_1(\Delta,\nu))=(1+o(1))\nu\times
\label{eq:fe}\end{equation}
$$\sum\limits_{\mu=0}^1\left(\frac{\Delta+(-1)^\mu}2
\log\left(\frac\Delta{\Delta+(-1)^\mu}\right)+
(-1)^\mu\frac\pi2\sum\limits_{\kappa=1}^{\left[\frac{\Delta+(-1)^\mu}2\right]}
\cot\left(\frac{\pi\kappa}{\Delta+(-1)^\mu}\right)\right).$$
In view of (\ref{eq:fd}),
\begin{equation}
d^\ast_0(5,\nu)=\exp(6\nu(\Delta+1)+O(\nu/\log(\nu))),
d^\ast_0(7,\nu)=
\label{eq:ff}\end{equation}
$$\exp(8\nu(8)+O(\nu/\log(\nu))).$$
In view of (\ref{eq:ff})
\begin{equation}
\log(d^\ast_1(5,\nu))=(1+o(1))\nu\times
\label{eq:fg}\end{equation}
$$(-3\log(1.2)+2\log(0.8)+
(\pi/2)(\cot(\pi/6)+\cot(\pi/3)+
\cot(\pi/4)))=$$
$$(1+o(1))\nu\times1.956124...,$$
\begin{equation}
\log(d^\ast_1(7,\nu))=(1+o(1))\nu\times
\label{eq:fh}\end{equation}
$$(4\log(7/8)+3\log(7/6))+
(1+o(1))(\pi/2)\nu\times$$
$$(-\cot(\pi/6)-\cot(\pi/3)+\cot(\pi/8)+\cot(3\pi/8)+\cot(\pi/4))=$$
$$(1+o(1))\nu(4\log(7/8)+3\log(7/6)+\pi(-2/\sqrt{3}+2/\sqrt{2}+1/2)=$$
$$(1+o(1))\nu\times2.314407\ldots\,,$$
when $\nu\to\infty.$

In view of (\ref{eq:ab}) and (\ref{eq:ac}),
$$\alpha^\ast(z;\nu)d^\ast_0(\nu)/d^\ast_1(\nu)\in{\mathbb Z}[z],$$
$$\phi^\ast(z;\nu)d^\ast_0(\nu)/d^\ast_1(\nu)\in{\mathbb Z}[z].$$
Let
\begin{equation}
U_\Delta(m,\nu)=d^\ast_0(\nu)/d^\ast_1(\nu),\,\Lambda_0(m)=0,
\label{eq:fi}\end{equation}
if
$m\ne2p^\alpha,$ where $p$ run over the all the prime numbers and $\alpha$
run over ${\mathbb N}$ and let
\begin{equation}\label{eq:gj}
U_\Delta(m,\nu)=\frac{d^\ast_0(\nu)}{d^\ast_1(\nu)}p^{[(\Delta+1)\nu/\phi(m)]+1}
,\,\Lambda_0(m)=\Lambda(m/2),
\end{equation}
if
$m=2p^\alpha,$ where $p$ is a prime number and $\alpha\in\mathbb N.$
In view of the (\ref{eq:ab}), (\ref{eq:ac}) and Lemma 1,
\begin{equation}\label{eq:gj1}
\alpha^\ast(z;\nu)\bigg\vert_
{z=\left(\frac1{2\cos(\frac{k\pi i}m)},\,\frac{k\pi i}m-\pi\right)}
U_\Delta(m,\nu)
\in{\mathbb Z}_{{\mathbb Q}[\exp(2i\pi/m)]},
\end{equation}
\begin{equation}\label{eq:gj2}
\phi^\ast(z;\nu)\bigg\vert_
{z=\left(\frac1{2\cos(\frac{k\pi i}m)},\,\frac{k\pi i}m-\pi\right)}
U_\Delta(m,\nu)
\in{\mathbb Z}_{{\mathbb Q}[\exp(2i\pi/m)]},
\end{equation}
 where $(k,m)=1.$
In view of (\ref{eq:gj}), (\ref{eq:fi}), (\ref{eq:fe}), (\ref{eq:fd}),
(\ref{eq:1}) and (\ref{eq:3})
\begin{equation}
\frac{d^\ast_0(\nu)}{d^\ast_1(\nu)}=
\label{eq:gb}\end{equation}
$$\nu(1+o(1))V^\ast_\Delta
\log(U_\Delta(m,\nu))=\nu(1+o(1))V_\Delta(m),$$
when $\nu\to\infty.$

The polynomial (\ref{eq:bj1}) take the form
$$
D^\wedge (z,\eta)=(\eta+1)\left(\eta+\frac{\Delta-1}{\Delta+1}\right)+
\frac{2\Delta\exp(i\psi)\eta}{(\Delta+1)\cos(\psi)}=
$$
$$((\Delta+1)\eta^2+2\Delta(2+iT)\eta+(\Delta-1))/(\Delta+1),$$
where $\psi\in(-pi/2,\,\pi/2)$ and $T=\tan(\psi);$
its roots are equal to
\begin{equation}
-(2\Delta+\Delta iT+R)/(\Delta+1),
\label{eq:hf}\end{equation}
where
$R^2=\Delta^2(3-T^2)+1+4\Delta^2iT.$
In view of (\ref{eq:0}),
Then
$$R\in\{\pm\left( w_\Delta(T)+i2\Delta^2iT/w_\Delta(T)\}\right).$$
In view of (\ref{eq:hf}) and (\ref{eq:ch}),
$$
 \eta_j^ \wedge (r,\psi,\delta_0)=
$$
$$
-\frac{
2\Delta+\Delta iT+(-1)^j\left( w_\Delta(T)+i2\Delta^2iT/w_\Delta(T)\right)}
{\Delta+1}=
$$
$$
-\frac{
2\Delta+(-1)^jw_\Delta(T)+iT\Delta\left(1+(-1)^j2\Delta/w_\Delta(T)\right)}
{\Delta+1},
$$
where $j=0,1,$
$$
\vert\eta_j^ \wedge (r,\psi,\delta_0)+k\vert^2=
$$
$$
\frac{\left(2\Delta+(-1)^jw_\Delta(T)-k(\Delta+1)\right)^2+
T^2\Delta^2\left(1+(-1)^j2\Delta/w_\Delta(T)\right)^2}
{(\Delta+1)^2},
$$
where $j=0,1;\,k=0,\,1,\,-1.$
Therefore, in view of (\ref{eq:ba}) and (\ref{eq:2})
\begin{equation}
\ln\vert h^\sim(\eta_j^ \wedge (r,\psi,\delta_0))\vert=
\label{eq:hg1}\end{equation}
$$(\eta_j(r,\psi,\delta_0)-1)(1-\delta_0)^{-d_1}
(\eta_j(r,\psi,\delta_0)+1)2^{-2}\eta_j(r,\psi,\delta_0)^{d_1}=$$
$$-log\left(4(\Delta+1)^{\Delta+1}(1-1/\Delta)^(\Delta-1)\right)+$$
$$\frac12\log\left(\left(2\Delta+(-1)^jw_\Delta(T)+(\Delta+1)\right)^2+
T^2\Delta^2\left(1+\frac{(-1)^j2\Delta}{w_\Delta(T)}\right)^2\right)+$$
$$\frac12
\log\left(\left(2\Delta+(-1)^jw_\Delta(T)-(\Delta+1)\right)^2+
T^2\Delta^2\left(1+\frac{(-1)^j2\Delta}{w_\Delta(T)}\right)^2\right)+$$
$$\frac{(\Delta-1)}2
\log\left(\left(2\Delta+(-1)^jw_\Delta(T)\right)^2+
T^2\Delta^2\left(1+\frac{(-1)^j2\Delta}{w_\Delta(T)}\right)^2\right)=$$
$$l_\Delta(j,T),$$
where $j=0,1.$
Clearly,
$$w_\Delta(0)=\sqrt{3\Delta^2+1},$$
$$
 \eta_j^ \wedge (1/2,0,\delta_0)=-\frac{2\Delta+(-1)^j\sqrt{3\Delta^2+1}}
{\Delta+1},
$$
where $j=0,1,$
$$
\left\vert\eta_j^ \wedge (1/2,0,\delta_0)+k\right\vert=
\left\vert\frac{2\Delta+(-1)^j\sqrt{3\Delta^2+1}-k(\Delta+1)}
{\Delta+1}\right\vert,
$$
where $j=0,1;\,k=0,\,1,\,-1.$
Therefore
\begin{equation}
l_\Delta(\epsilon,0)=
\left(\log\vert h^\sim(\eta_\epsilon^ \wedge(1/2,0,\delta_0))\vert
\right)=\label{eq:h11}\end{equation}
$$
\log\left(\vert(\eta_\epsilon(1/2,0,\delta_0)-1)(1-\delta_0)^{-d_1}
(\eta_\epsilon(1/2,0,\delta_0)+1)
2^{-2}\eta_\epsilon(1/2,0,\delta_0)^{d_1}\vert\right)=
$$
$$-\log\left(4(\Delta+1)^{\Delta+1}(1-1/\Delta)^(\Delta-1)\right)\,+$$
$$\log\left(\vert2\Delta+(-1)^\epsilon\sqrt{3\Delta^2+1}-(\Delta+1)\vert
\right)+$$
$$\log\left(\vert{2\Delta+(-1)^\epsilon\sqrt{3\Delta^2+1}+(\Delta+1)}\vert
\right)+$$
$$(\Delta-1)\log\left(\vert{2\Delta+(-1)^\epsilon\sqrt{3\Delta^2+1}}\vert
\right).$$
Consequently
$$l_5(1,0)=-\|og(4)-6\log6-4\log(0.8)+$$
$$
\log(\sqrt{76}-4)+\log(16-\sqrt{76})+4\log(10-\sqrt{76})$$
I made computations below "by hands" using calculator of the firm "CASIO."
$$\log4=1,386294361...\,;\,6\log(6)=10,7505682...\,;$$
$$4\log(0.8)=-0,892574205...\,;$$
$$\sqrt{76}=8,717797887...\,;\,\sqrt{76}-4=4,717797887...\,;$$
$$16-\sqrt{76}=7,282202113...\,;\,10-\sqrt{76}=1,282202113...\,;$$
$$\log\left(\sqrt{76}-4\right)=1.551342141...\,;\,
\log\left(16-\sqrt{76}\right)=1.985433305...\,;$$
$$\log\left(10-\sqrt{76}\right)=0.248579...\,;\,
4\log\left(10-\sqrt{76}\right)=0,994316001...\,;\,
$$
\begin{equation}
l_5(1,0)=-6.713196909...;
\label{eq:h12}\end{equation}
$$
l_7(1,0)=-\log(4)-8\log(8)-6\log(6)+6\log(7)+$$
$$\log\left(\sqrt{148}-6\right)+
\log\left(22-\sqrt{148}\right)+
6\log\left(14-\sqrt{148}\right);$$
$$8\log8=16,63553233...\,;\,6\log6=10,75055682...;\,6\log7=11,67546089...;$$
$$\sqrt{148}=12,16552506...\,;\,\sqrt{148}-6=6,16552506...\,$$
$$22-\sqrt{148}=9,83474939...\,;\,14-\sqrt{148}=1,83474939...\,;$$
$$\log(\sqrt{148}-6)=1,818973301;\,\log(22-\sqrt{148})=2,285894063...;\,$$
$$\log(14-\sqrt{148})=0,606758304...\,;\,6\log(14-\sqrt{148})=3,640549824...\,;
\,$$
\begin{equation}
l_7(1,0)=-9,35150543...\,.
\label{eq:h13}\end{equation}
In view of (\ref{eq:1}),\,(\ref{eq:fd}),\,(\ref{eq:fe}),\,(\ref{eq:fg}),\,
(\ref{eq:fh}) and (\ref{eq:gb}),
\begin{equation}
V_5^\ast=6-1.956124...=4,04387...;V_7^\ast=8-2.314407=5,685593.
\label{eq:h14}\end{equation}

In view (\ref{eq:h12}) -- (\ref{eq:h14}),
\begin{equation}
-V_5^\ast-l_5(1,0)>0,\,-V_7^\ast-l_7(1,0)>0.
\label{eq:h15}\end{equation}
So, the key inequalities (\ref{eq:h15}) are checked "by hands".
I view of (\ref{eq:hg1}), (\ref{eq:h15}) and Lemma 3,
$$
-V_5^\ast-l_5(1,\tan(\pi/m))>0,\,-V_7^\ast-l_7(1,\tan(\pi/m))>0,
$$
where $m>2.$ Since $(\log(p))/(p^{\alpha-1}(p-1))$ decreases together with
increasing of $p\in(3,\,+\infty)$
with fixed $\alpha\ge1,$ or icreasing of $\alpha\in(1,\,+\infty$
with fixed $p\ge2$ (or, of course, increasing both $\alpha\in(1,\,+\infty$
and $p\in(3,\,+\infty)$), and
$$\lim\limits_{p\to\infty}((\log(p))/(p^{\alpha-1}(p-1)))=0,$$
where $\alpha\ge1,$
$$\lim\limits_{\alpha\to\infty}((\log(p))/(p^{\alpha-1}(p-1)))=0,$$
where $p\ge2,$ it follows that the inequality (\ref{eq:4}) holds for all the
sufficient big integers $m.$ Computations on computer of class "Pentium" show
that the inequality (\ref{eq:4}) holds for $m=3,\,m=4,\,m=5$
and $m=2\times5;$ therefore inequality (\ref{eq:4})
 holds for all the $m>2\times3.$
Let
$
\varepsilon_0=h_{\Delta}(m)/2,
$
with $h_{\Delta}(m)$ defined in (\ref{eq:2a}).
 In  view of (\ref{eq:4}),
$
\varepsilon_0>0.
$
We take now $K=K_m={\mathbb Q}[\exp(2\pi i/m)].$
Let further $\{\sigma_1,\,\ldots,\,\sigma_{\phi(m)}\}=\Gal(K/{\mathbb Q}).$ For each $j=1,\,\ldots,\,\phi(m)$
there exists $k_j\in(-m/2,m/2)\cap{\mathbb Z}$ such that
$$(\vert k_j\vert,\,m)=1,\,
\sigma_j\left(\exp\left(\frac{2\pi i}m\right)\right)=
\exp\left(\frac{2\pi ik_j}m\right).$$
Let $a$ be the element of $\mathfrak K,$ such that
$$
\pi_j(a)=\log(2+\sigma_j(\exp(2\pi i/m)))=\log(2+\exp(2\pi ik_j/m)),
$$
where $j=1,\,\ldots,\,\phi(m);$ we suppose that $k_1=1.$
In view of (\ref{eq:gj1}) and (\ref{eq:gj2}),
let $\alpha_1^\vee (\nu),\,\alpha_1^\wedge (\nu),\,
\alpha_2^\vee (\nu),\,\alpha_2^\wedge (\nu),\,$
 are elements
in ${\mathfrak K}$ such that
$$
\pi_j(\alpha_1^\vee(\nu))=
\alpha^\ast(z;\nu)\bigg\vert_
{z=\left(\frac1{2\cos(\frac{k_j\pi i}m)},\,\frac{k_j\pi i}m-\pi\right)},
$$
$$
\pi_j(\alpha_2^\vee(\nu))=
\phi^\ast(z;\nu)\bigg\vert_
{z=\left(\frac1{2\cos(\frac{k_j\pi i}m)},\,\frac{k_j\pi i}m-\pi\right)},
$$
\begin{equation}\label{eq:ij}
\pi_j(\alpha_1^\wedge(\nu))=
\alpha^\ast(z;\nu)\bigg\vert_
{z=\left(\frac1{2\cos(\frac{k_j\pi i}m)},\,\frac{k_j\pi i}m-\pi\right)}
U_\Delta(m,\nu),
\end{equation}
\begin{equation}\label{eq:ia}
\pi_j(\alpha_2^\wedge(\nu))=
\phi^\ast(z;\nu)\bigg\vert_
{z=\left(\frac1{2\cos(\frac{k_j\pi i}m)},\,\frac{k_j\pi i}m-\pi\right)}
U_\Delta(m,\nu),
\end{equation}
where $j=1,\,\ldots,\,\phi(m).$
Then $\alpha_k^\wedge (\nu)\in{\mathbb Z}_K$ for $k=1,\,2.$
\begin{equation}\label{eq:ec00}
y^\vee(\nu)=-\alpha^\vee_2(\nu)+a\alpha_1^\vee(\nu),
\end{equation}
and let $y(\nu)$ is defined by means the equality (\ref{eq:ec}).
According to the Corrollary of the Lemma 4, to the Theorem 4 in [\ref{r:cb}]
 (or Theorem 7 in [\ref{r:cf1}]),
 to the Lemma 8, to (\ref{eq:hg1}),
there exist $m^\ast_1\in{\mathbb N}$ having the following property:

for any $\varepsilon\in(0,\,\varepsilon_0)$ there exist
 $\gamma_0(\varepsilon)>0,\,\gamma_1(\varepsilon)>0,$  
and $\gamma_2(\varepsilon)>0$ such that
\begin{equation}\label{eq:ib0}
\vert\pi_j(\alpha_k^\vee(\nu))\vert\le
\end{equation}
$$
\gamma_0(\varepsilon)\exp((l_\Delta(\tan((k_j\pi i)/m),0)+
\varepsilon/3)\nu),
$$
where $k=1,\,2,\,j=1,\,\ldots,\,\phi(m)$ and $\nu\in{\mathbb N}-1+m^\ast_1,$
 \begin{equation}\label{eq:ib}
\gamma_1(\varepsilon)\exp((l_\Delta(\tan((k_j\pi i)/m),1)-
\varepsilon/3)\nu)\le
\end{equation}
$$\max(\vert\pi_j(y^\vee(\nu))\vert,\,\vert\pi_j(y^\vee(\nu+1))\vert\le$$
$$\gamma_2(\varepsilon)\exp((l_\Delta(\tan((k_j\pi i)/m),1)+
\varepsilon/3)\nu),$$
where $j=1,\,\ldots,\,\phi(m)$ and $\nu\in{\mathbb N}-1+m^\ast_1.$

Let $\omega_1(m)=(m-1)/2,$ if $m$ is odd,
$\omega_1(m)=m/2-2,$ if $m\equiv2(\mod4)$ and $\omega(m)=m/2-1,$
if $m\equiv0(\mod4).$ Then
 $$\omega_1(m)=\sup\{k\in{\mathbb N}\colon k_j<m/2,(k,m)=1\}.$$

According to the Lemma 3 and (\ref{eq:hg1}),
\begin{equation}\label{eq:ic0}
l_\Delta(\tan((k_j\pi i)/m),0)\le l_\Delta(\tan((\omega_1(m)\pi i)/m),0),
\end{equation}
\begin{equation}\label{eq:ic}
l_\Delta(\tan((\omega_1(m)\pi i)/m),1)\le\end{equation}
$$l_\Delta(\tan((k_j\pi i)/m),1)\le l_\Delta(\tan((\pi i)/m),1)$$
where $j=1,\,\ldots,\,\phi(m).$
In view of (\ref{eq:ib0}) -- (\ref{eq:ic}),
\begin{equation}\label{eq:id0}
\vert\pi_j(\alpha_k^\vee(\nu))\vert\le
\gamma_0(\varepsilon)\exp((l_\Delta(\tan((\omega_1(\nu)\pi i)/m),0)+
\varepsilon/3)\nu),\end{equation}
where $k=1,\,2,\,j=1,\,\ldots,\,\phi(m)$ and $\nu\in{\mathbb N}-1+m^\ast_1,$
 \begin{equation}\label{eq:id}
\gamma_1(\varepsilon)\exp((l_\Delta(\tan((\omega_1(m)\pi i)/m),1)-
\varepsilon/3)\nu)\le
\end{equation}
$$\max(\vert\pi_j(y^\vee(\nu))\vert,\,\vert\pi_j(y^\vee(\nu+1))\vert\le$$
$$\gamma_2(\varepsilon)\exp((l_\Delta(\tan((\pi i)/m),1)+
\varepsilon/3)\nu),$$
where $j=1,\,\ldots,\,\phi(m)$ and $\nu\in{\mathbb N}-1+m^\ast_1.$
In view of (\ref{eq:gb}), there exists 
$m^\ast_2\in{\mathbb N}-1+m^\ast_1,$ such that
\begin{equation}\label{eq:ie}
\exp(V_\Delta(m)-\varepsilon/3)\nu\le U_\Delta(m,\nu)\le
\exp(V_\Delta(m)-\varepsilon/3)\nu
\end{equation}
where $\nu\in{\mathbb N}-1+m^\ast_2.$

In view of (\ref{eq:ic}) -- (\ref{eq:ie}), (\ref{eq:ij}) -- (\ref{eq:ec00}),
  (\ref{eq:2a}), (\ref{eq:2a0}),
\begin{equation}\label{eq:id1}
\vert\pi_j(\alpha_k(\nu))\vert\le
\gamma_0(\varepsilon)\exp((g_{\Delta,0}(m)+2\varepsilon/3)\nu),
\end{equation}
where $k=1,\,2,\,j=1,\,\ldots,\,\phi(m)$ and $\nu\in{\mathbb N}-1+m^\ast_2,$
 \begin{equation}\label{eq:id2}
\gamma_1(\varepsilon)\exp((-g_{\Delta,1}(m)-2\varepsilon/3)\nu)\le
\end{equation}
$$\max(\vert\pi_j(y^\vee(\nu))\vert,\,\vert\pi_j(y^\vee(\nu+1))\vert\le$$
$$\gamma_2(\varepsilon)\exp((-h_{\Delta}(m)+2\varepsilon/3)\nu),$$
where $j=1,\,\ldots,\,\phi(m)$ and $\nu\in{\mathbb N}-1+m^\ast_2.$

Let $X\in{\mathbb Z}_{K_m}\diagdown\{0\}.$ Then, in view of
(\ref{eq:id1}) and (\ref{eq:id2}),
\begin{equation}\label{eq:id3}
\vert\pi_j(X\alpha_k(\nu))\vert\vert\le
\gamma_0(\varepsilon)\exp((g_{\Delta,0}(m)+2\varepsilon/3)\nu)
\vert\pi_j(X)\vert\le
\end{equation}
$$
\gamma_0(\varepsilon)\exp((g_{\Delta,0}(m)+2\varepsilon/3)\nu)
q_\infty^{({\mathfrak K)}}(X),
$$
where $k=1,\,2,\,j=1,\,\ldots,\,\phi(m)$ and $\nu\in{\mathbb N}-1+m^\ast_2,$
 \begin{equation}\label{eq:id4}
\gamma_1(\varepsilon)\exp((-g_{\Delta,1}(m)-2\varepsilon/3)\nu)
\vert\pi_j(X)\vert\le
\end{equation}
$$\max(\vert\pi_j(Xy^\vee(\nu))\vert,\,\vert\pi_j(Xy^\vee(\nu+1))\vert\le$$
$$\max(q_\infty^{({\mathfrak K})}(Xy^\vee(\nu)),\,
q_\infty^{({\mathfrak K})}(Xy^\vee(\nu+1)),$$
where $j=1,\,\ldots,\,\phi(m)$ and $\nu\in{\mathbb N}-1+m^\ast_2,$
\begin{equation}\label{eq:id5}
\max(\vert\pi_j(Xy^\vee(\nu))\vert,\,\vert\pi_j(Xy^\vee(\nu+1))\vert\le
\end{equation}
$$\gamma_2(\varepsilon)\exp((-h_{\Delta}(m)+2\varepsilon/3)\nu)\vert\pi_j(X)
\vert\le$$
$$\gamma_2(\varepsilon)\exp((-h_{\Delta}(m)+2\varepsilon/3)\nu)
q_\infty^{\mathfrak K}(X),
$$
where $j=1,\,\ldots,\,\phi(m)$ and $\nu\in{\mathbb N}-1+m^\ast_2.$

In view of (\ref{eq:id3})
\begin{equation}\label{eq:id6}
q_\infty^{({\mathfrak K})}(X\alpha_k(\nu))\le
\end{equation}
$$
\gamma_0(\varepsilon)\exp((g_{\Delta,0}(m)+2\varepsilon/3)\nu)
q_\infty^{({\mathfrak K})}(X),
$$
where $k=1,\,2,\,$ and $\nu\in{\mathbb N}-1+m^\ast_2.$
In view of (\ref{eq:id5}),
 \begin{equation}\label{eq:id7}
\max(q_\infty^{({\mathfrak K})}(Xy^\vee(\nu)),\,
q_\infty^{({\mathfrak K})}((Xy^\vee(\nu+1))=
\end{equation}
$$
\sup(\{\vert\pi_j(Xy^\vee(\nu+\epsilon))\vert
,\,\colon \epsilon\in\{0,\,1\},\,j=1,\,\ldots,\,\phi(m)\})\le$$
$$\gamma_2(\varepsilon)\exp((-h_{\Delta}(m)+2\varepsilon/3)\nu)
q_\infty^{({\mathfrak K})}(X),
$$
where $\nu\in{\mathbb N}-1+m^\ast_2.$

Taking in acount (\ref{eq:id6}), (\ref{eq:id7}) and (\ref{eq:id4}),
we see that all the conditions of the Corollary of the Lemma 12
are fulfilled for 
$$\varepsilon\in(0,\,\varepsilon_0),
\gamma_0(\varepsilon),\,\gamma_1(\varepsilon),\,
\gamma_2(\varepsilon),
y=y(\nu),\,\alpha_1(\nu),\alpha_2(\nu),$$
$$r_1=r_1(\varepsilon)=\exp(g_{\Delta,0}(m)+2\varepsilon/3,$$
$$R_1=R_1(\varepsilon)=\exp(g_{\Delta,1}(m)+2\varepsilon/3),$$
$$R_2=R_2(\varepsilon)\exp(h_{\Delta}(m)-2\varepsilon/3),$$
and this proves the part of our Theorem connected with
the inequality (\ref{eq:5}).

Let again $X\in{\mathbb Z}_{K_m}\diagdown\{0\}$ and let
$$
q_{min}^{({\mathfrak K})}(X)=\inf(\vert\{\pi_j(X)\vert
\colon j=1,\,\ldots,\,\phi(m)\})
$$
Clearly, $q_{min}^{({\mathfrak K})}(X)>0$ According to the Theorem 4
 in [\ref{r:cb}], or to the Theorem 7 in [\ref{r:cf1}],
there exist $m^\ast_1\in{\mathbb N}$ having the following property:
for any $\varepsilon\in(0,\,\varepsilon_0)$ there exist
 $\gamma_0^\ast(X,\varepsilon)>0,\,\gamma_1^\ast(X,\varepsilon)>0,$  
and $\gamma_2^\ast(X,\varepsilon)>0$ such that
$$
\vert\pi_j(\alpha_k^\vee(\nu))\vert\le\gamma_0^\ast
(\varepsilon)\exp((l_\Delta(\tan((\omega_m\pi i)/m),0)+\varepsilon/3)\nu),
$$
where $k=1,\,2,\,j=1,\,\ldots,\,\phi(m)$ and $\nu\in{\mathbb N}-1+m^\ast_1,$
$$
\gamma_1^\ast(X\varepsilon)\exp((l_\Delta(\tan((\pi i)/m),1)-
\varepsilon/3)\nu)\le
$$
$$\max(\vert\pi_j(y^\vee(\nu))\vert,\,\vert\pi_j(y^\vee(\nu+1))\vert\le$$
$$\gamma_2(\varepsilon)\exp((l_\Delta(\tan((\pi i)/m),1)+
\varepsilon/3)\nu),$$
where $j=1,\,\ldots,\,\phi(m)$ and $\nu\in{\mathbb N}-1+m^\ast_1.$
Repeating the previous considerations, we see
that all the conditions of the Corollary of the Lemma 12
are fulfilled for $\varepsilon\in(0,\,\varepsilon_0),$
$$\gamma_0=\gamma_0^\ast(X,\varepsilon),\,
\gamma_1=\gamma_1^\ast(X,\varepsilon),\,
\gamma_2=\gamma_2^\ast(X,\varepsilon),$$
$$y=y(\nu),\,\alpha_1(\nu),\alpha_2(\nu),\,
r_1=r_1(\varepsilon)=\exp(g_{\Delta,0}(m)+2\varepsilon/3,$$
and
$$R_1=R_2=R_2(\varepsilon)\exp(h_{\Delta}(m)-2\varepsilon/3),$$
and this proves the part of our Theorem connected with
the inequality (\ref{eq:6}). $\blacksquare$

Below are values of $\beta$ and $\alpha$ computed
for $\Delta\in\{5,\,7\}$ and some $m\in{\mathbb N}.$
$$(m;\,\Delta;\,\beta;\,\alpha)=(3;\,5;\,3,111228...\,;\,3,111228...),$$
$$(m;\,\Delta;\,\beta;\,\alpha)=(3;\,7;\,3,073525...\,;\,3,073525...),$$
$$(m;\,\Delta;\,\beta;\,\alpha)=(4;\,5;\,11,458947...\,;\,11,458947...),$$
$$(m;\,\Delta;\,\beta;\,\alpha)=(4;\,7;\,10,551730...\,;\,10,551730...),$$
$$(m;\,\Delta;\,\beta;\,\alpha)=(5;\,5;\,4,826751...\,;\,5,607961...),$$
$$(m;\,\Delta;\,\beta;\,\alpha)=(5;\,7,\,4,837858...\,;\,5,684622...),$$
$$(m;\,\Delta;\,\beta;\,\alpha)=(7;\,5;\,5,701485...\,;\,6,977258...),$$
$$(m;\,\Delta;\,\beta;\,\alpha)=(7;\,7;\,5,724804...\,;\,7,114963...),$$
$$(m;\,\Delta;\,\beta;\,\alpha)=(8;\,5;\,8,337857...\,;\,9,436901...),$$
$$(m;\,\Delta;\,\beta;\,\alpha)=(8;\,7;\,8,253047...\,;\,9,433260...),$$
$$(m;\,\Delta;\,\beta;\,\alpha)=(9;\,5;\,6,312056...\,;\,7,960502...),$$
$$(m;\,\Delta;\,\beta;\,\alpha)=(9;\,7;\,6,335274...\,;\,8,134962...),$$
$$(m;\,\Delta;\,\beta;\,\alpha)=(10;\,5;\,43,546644...\,;\,46,230614...),$$
$$(m;\,\Delta;\,\beta;\,\alpha)=(10;\,7;\,35,648681...\,;\,38,043440...),$$
$$(m;\,\Delta;\,\beta;\,\alpha)=(11;\,5;\,6,786990...\,;\,8,735234...),$$
$$(m;\,\Delta;\,\beta;\,\alpha)=(11;\,7,\,6,806087...\,;\,8,934922...),$$
$$(m;\,\Delta;\,\beta;\,\alpha)=(12;\,5;\,5,638541...\,;\,6,813222...),$$
$$(m;\,\Delta;\,\beta;\,\alpha)=(12;\,7;\,5,696732...\,;\,6,983870...),$$
$$(m;\,\Delta;\,\beta;\,\alpha)=(13;\,5;\,7,177155...\,;\,9,376030...),$$
$$(m;\,\Delta;\,\beta;\,\alpha)=(13;\,7;\,7,190814...\,;\,9,594580...),$$
$$(m;\,\Delta;\,\beta;\,\alpha)=(14;\,5;\,19,659885...\,;\,21,835056...),$$
$$(m;\,\Delta;\,\beta;\,\alpha)=(14;\,7;\,18,447228...\,;\,20,668254...),$$
$$(m;\,\Delta;\,\beta;\,\alpha)=(15\,;5;\,7,508714...\,;\,9,922761...),$$
$$(m;\,\Delta;\,\beta;\,\alpha)=(15;\,7;\,7,516606...\,;\,10,156245...),$$
$$(m;\,\Delta;\,\beta;\,\alpha)=(16;\,5,\,7,951153...\,;\,9,876454...),$$
$$(m;\,\Delta;\,\beta;\,\alpha)=(16;\,7,\,7,945763...\,;\,10,039605...),$$
$$(m;\,\Delta;\,\beta;\,\alpha)=(17;\,5;\,7,797153...\,;\,10,399610...),$$
$$(m;\,\Delta;\,\beta;\,\alpha)=(17;\,7,\,7,799343...\,;\,10,645404...),$$
$$(m;\,\Delta;\,\beta;\,\alpha)=(18;\,5,\,9,486110...\,;\,10,955534...),$$
$$(m;\,\Delta;\,\beta;\,\alpha)=(18;\,7,\,9,406368...\,;\,10,989150...),$$
$$(m;\,\Delta;\,\beta;\,\alpha)=(19;\,5;\,8,052478...\,;\,10,822446...),$$ 
$$(m;\,\Delta;\,\beta;\,\alpha)=(19;\,7;\,8,049182...\,;\,11,078690...),$$
$$(m;\,\Delta;\,\beta;\,\alpha)=(20;\,5;\,6,696241...\,;\,8,559091...),$$
$$(m;\,\Delta;\,\beta;\,\alpha)=(20;\,7;\,6,733979...\,;\,8,774063...),$$
$$(m;\,\Delta;\,\beta;\,\alpha)=(21;\,5;\,8,281548...\,;\,11,202268...),$$
$$(m;\,\Delta;\,\beta;\,\alpha)=(21;\,7;\,8,273039...\,;\,11,467583...),$$ 
$$(m;\,\Delta;\,\beta;\,\alpha)=(22;\,5;\,13,134623...\,;\,15,504916...),$$
$$(m;\,\Delta;\,\beta;\,\alpha)=(22;\,7;\,12,815391...\,;\,15,331975...),$$
$$(m;\,\Delta;\,\beta;\,\alpha)=(23;\,5;\,8,489281...\,;\,11,547024...),$$
$$(m;\,\Delta;\,\beta;\,\alpha)=(23;\,7;\,8,475843...\,;\,11,820351...),$$
$$(m;\,\Delta;\,\beta;\,\alpha)=(24;\,5;\,7,088338...\,;\,9,210037...),$$
$$(m;\,\Delta;\,\beta;\,\alpha)=(24;\,7;\,7,116679...\,;\,8,439782...),$$
$$(m;\,\Delta;\,\beta;\,\alpha)=(25;\,5;\,8,679328...\,;\,11,862643...),$$
$$(m;\,\Delta;\,\beta;\,\alpha)=(25;\,7;\,8,661235...\,;\,12,143143...),$$
$$(m;\,\Delta;\,\beta;\,\alpha)=(26;\,5;\,12,172520...\,;\,14,674949...),$$
$$(m;\,\Delta;\,\beta;\,\alpha)=(26;\,7;\,11,944943...\,;\,14,618461...),$$
$$\ldots$$
$$(m;\,\Delta;\,\beta;\,\alpha)=(32;\,5;\,8,654733...\,;\,11,466214...),$$
$$(m;\,\Delta;\,\beta;\,\alpha)=(32;\,7;\,8,637697...\,;\,11,705492...),$$
$$(m;\,\Delta;\,\beta;\,\alpha)=(33;\,5;\,9,310125...\,;\,12,911341...),$$
$$(m;\,\Delta;\,\beta;\,\alpha)=(33;\,5;\,9,275806...\,;\,13,214792...),$$
{\begin{center}\large\bf References.\end{center}}
\footnotesize
\vskip4pt
\refstepcounter{r}\noindent[\ther]
R.Ap\'ery, Interpolation des fractions continues\\
\hspace*{3cm}  et irrationalite de certaines constantes,\\
\hspace*{3cm} Bulletin de la section des sciences du C.T.H., 1981, No 3,
 37 -- 53;
\label{r:cd}\\
\refstepcounter{r}
\noindent[\ther]
F.Beukers, A note on the irrationality of $\zeta (2)$ and $\zeta (3),$\\
\hspace*{3cm} Bull. London Math. Soc., 1979, 11, 268 --  272;
\label{r:ce}\\
\refstepcounter{r}
\noindent[\ther]
 A.van der Porten,
 A proof that Euler missed...Ap\'ery's proof of the irrationality of
 $\zeta (3),$\\
\hspace*{3cm} Math Intellegencer, 1979, 1, 195 -- 203;\label{r:cf}\\
\refstepcounter{r}
\noindent[\ther]
 W. Maier, Potenzreihen irrationalen Grenzwertes,\\
\hspace*{3cm} J.reine angew. Math.,
 156, 1927, 93 -- 148;\label{r:cg}\\
\refstepcounter{r}\noindent[\ther]
 E.M. Niki\u sin, On irrationality of the values of the functions F(x,s)
 (in Russian),\\
\hspace*{3cm} Mat.Sb. 109 (1979), 410 -- 417;\\
\hspace*{3cm}
 English transl. in Math. USSR Sb. 37 (1980), 381 -- 388;\label{r:ch}\\
\refstepcounter{r}\noindent[\ther]
 G.V. Chudnovsky, Pade approximations to the generalized
 hyper-geometric functions\\
\hspace*{3.4cm}  I,J.Math.Pures Appl., 58, 1979,
 445 -- 476;\label{r:dj}\\
\refstepcounter{r}\noindent[\ther]
\rule{2.7cm}{.3pt},
 Transcendental numbers,  Number Theory,Carbondale,\\
\hspace*{3.4cm}  Lecture Notes in Math, Springer-Verlag, 1979, 751, 45 -- 69;
\label{r:da}\\
\refstepcounter{r}\noindent[\ther]
\rule{2.7cm}{.3pt}, Approximations rationelles des logarithmes
 de nombres rationelles\\
\hspace*{3.4cm} C.R.Acad.Sc. Paris, S\'erie A, 1979, 288, 607 -- 609;
\label{r:db}\\
\refstepcounter{r}\noindent[\ther]
\rule{2.7cm}{.3pt}, Formules d'Hermite pour les approximants de Pad\'e de
logarithmes\\
\hspace*{3.4cm} et de fonctions bin\^omes, et mesures d'irrationalit\'e,\\
\hspace*{3.4cm} C.R.Acad.Sc. Paris, S\'erie A, 1979, t.288, 965 -- 967;
\label{r:dc}\\
\refstepcounter{r}\noindent[\ther]
\rule{2.5cm}{.3pt},Un syst\'me explicite d'approximants de Pad\'e\\
\hspace*{3.4cm} pour les fonctions hyp\'erg\'eometriques g\'en\'eralies\'ees,\\
\hspace*{3.4cm} avec applications a l'arithm\'etique,\\
\hspace*{3.4cm}
 C.R.Acad.Sc. Paris, S\'erie A, 1979, t.288, 1001 -- 1004;\label{r:dd}\\
\refstepcounter{r}\noindent[\ther]
\rule{2.5cm}{.3pt},
 Recurrenses defining Rational Approximations\\
\hspace*{3.4cm} to the irrational numbers, Proceedings\\
\hspace*{3.4cm}  of the Japan Academie,
 Ser. A, 1982, 58, 129 -- 133;
\label{r:de}\\
\refstepcounter{r}\noindent[\ther]
\rule{2.5cm}{.3pt}, On the method of Thue-Siegel,\\
\hspace*{3.4cm} Annals of Mathematics, 117 (1983), 325 -- 382;
\label{r:df}\\
\refstepcounter{r}\noindent[\ther]
K.Alladi  and M. Robinson, Legendre polinomials and irrationality,\\
\hspace*{6cm}J. Reine Angew.Math., 1980, 318, 137 -- 155;
\label{r:dg}\\
\refstepcounter{r}\noindent[\ther]
A. Dubitskas, An approximation of logarithms of some numbers,\\
\hspace*{3cm} Diophantine approximations II,Moscow, 1986, 20 -- 34;
\label{r:dh}\\
\refstepcounter{r}\noindent [\ther] \rule{2cm}{.3pt},
 On approximation of $\pi/ \sqrt {3}$ by rational fractions,\\
\hspace*{3cm} Vestnik MGU, series 1, 1987, 6, 73 -- 76;
\label{r:ej}\\
\refstepcounter{r}\noindent[\ther]
 S.Eckmann, \"Uber die lineare Unaqbhangigkeit der Werte gewisser Reihen,\\
\hspace*{3cm} Results in Mathematics, 11, 1987, 7 -- 43;
\label{r:ea}\\
\refstepcounter{r}\noindent[\ther]
 M.Hata, Legendre type polinomials and irrationality mesures,\\
\hspace*{3cm} J. Reine Angew. Math., 1990, 407, 99 -- 125;
\label{r:eb}\\
\refstepcounter{r}\noindent[\ther]
 A.O. Gelfond, Transcendental and algebraic numbers (in Russian),\\
\hspace*{3cm} GIT-TL, Moscow, 1952;
\label{r:ec}\\
\refstepcounter{r}\noindent[\ther]
H.Bateman and A.Erd\'elyi, Higher transcendental functions,1953,\\
\hspace*{3cm}  New-York -- Toronto -- London,
 Mc. Grow-Hill Book Company, Inc.;
\label{r:ed}\\
\refstepcounter{r}\noindent[\ther]
 E.C.Titchmarsh, The Theory of Functions, 1939, Oxford University Press;
\label{r:ee}\\
\refstepcounter{r}\noindent[\ther]
E.T.Whittaker and G.N. Watson, A course of modern analysis,\\
\hspace*{3cm} 1927, Cambridge University Press;
\label{r:ef}\\
\refstepcounter{r}\noindent[\ther]
 O.Perron, \"Uber die Poincaresche Differenzengleichumg,\\
\hspace*{3cm} Journal f\"ur die reine und angewandte mathematik,\\
\hspace*{3cm} 1910, 137,  6 -- 64;\label{r:a}\\
\refstepcounter{r}
\noindent [\ther] A.O.Gelfond, Differenzenrechnung (in Russian),
 1967, Nauka, Moscow.
\label{r:b}\\
\refstepcounter{r}\noindent [\ther]
 A.O.Gelfond and I.M.Kubenskaya, On the theorem of Perron\\
\hspace*{6cm} in the theory of differrence equations  (in Russian),\\
\hspace*{3cm} IAN USSR, math. ser., 1953, 17, 2,  83 --  86.
\label{r:c}\\
\refstepcounter{r}
\noindent [\ther] M.A.Evgrafov, New proof of the theorem of Perron\\
\hspace*{3cm} (in Russian),IAN USSR, math. ser., 1953, 17, 2, 77 --  82;
\label{r:d}\\
\refstepcounter{r}
\noindent [\ther] G.A.Frejman, On theorems of of Poincar\'e and Perron\\
\hspace*{3cm} (in Russian), UMN, 1957, 12, 3 (75), 243 -- 245;
\label{r:e}\\
\refstepcounter{r}
\noindent [\ther]  N.E.N\"orlund, Differenzenrechnung, Berlin,
 Springer Verlag, 1924;\label{r:f}\\
\refstepcounter{r}
\noindent [\ther] I.M.Vinogradov, Foundtions of the Number Theory,
(in Russian), 1952, GIT-TL;\label{r:g}\\
\refstepcounter{r}
\noindent [\ther] \rule{2cm}{.3pt},
J. Diedonne, Foundations of modern analysis,\\
\hspace*{3cm}
 Institut des Hautes \'Etudes Scientifiques, Paris,\\
\hspace*{3cm} Academic Press, New York and London, 1960\label{r:ei}\\
\refstepcounter{r}
\noindent [\ther] CH.-J. de la Vall\'ee Poussin, Course d'analyse
infinit\'esimale,\\
\hspace*{3cm} Russian translation by G.M.Fikhtengolts,\\
\hspace*{3cm} GT-TI, 1933;
\label{r:cb0}\\
\refstepcounter{r}
\noindent [\ther] \rule{2cm}{.3pt},
H.WEYL, Algebraic theory of numbers, 1940,\\
\hspace*{3cm} Russian translation by L.I.Kopejkina,
\label{r:cb1}\\
\refstepcounter{r}
\noindent [\ther] L.A.Gutnik, On the decomposition
 of the difference operators of Poincar\'e type\\
\hspace*{3cm} (in Russian), VINITI, Moscow, 1992, 2468 -- 92, 1 -- 55;
\label{r:h}\\
\refstepcounter{r}
\noindent [\ther] \rule{2cm}{.3pt},  On the decomposition
 of the difference operators\\
\hspace*{3cm}  of Poincar\'e type in Banach algebras\\
\hspace*{3cm} (in Russian), VINITI, Moscow, 1992, 3443 -- 92, 1 -- 36;
\label{r:i}\\
\refstepcounter{r}
\noindent [\ther] \rule{2cm}{.3pt}, On the difference equations of Poincar\'e
 type\\
\hspace*{3cm}  (in Russian), VINITI, Moscow 1993, 443 -- B93, 1 -- 41;
\label{r:aj}\\
\refstepcounter{r}
\noindent [\ther] \rule{2cm}{.3pt}, On the difference equations
 of Poincar\'e type in normed algebras\\
\hspace*{3cm} (in Russian), VINITI, Moscow, 1994, 668 -- B94, 1 -- 44;
\label{r:aa}\\
\refstepcounter{r}
\noindent [\ther] \rule{2cm}{.3pt}, On the decomposition of
 the difference equations  of Poincar\'e type\\
\hspace*{3cm} (in Russian),
 VINITI, Moscow, 1997, 2062 -- B97, 1 -- 41;
\label{r:ab}\\
\refstepcounter{r}
\noindent [\ther] \rule{2cm}{.3pt}, The difference equations
 of Poincar\'e type\\
\hspace*{3cm} with characteristic polynomial having roots equal to zero\\
\hspace*{3cm}(in Russian), VINITI, Moscow, 1997, 2418 -- 97, 1 -- 20;
\label{r:ac}\\
\refstepcounter{r}
\noindent [\ther]  \rule{2cm}{.3pt}, On the behavior of solutions\\
\hspace*{3cm} of difference equations of Poincar\'e type\\
\hspace*{3cm}
 (in Russian), VINITI, Moscow, 1997, 3384 -- B97, 1 - 41;
\label{r:ad}\\
\refstepcounter{r}
\noindent [\ther] \rule{2cm}{.3pt}, On the variability of solutions
of difference equations of Poincar\'e type\\
\hspace*{3cm}
 (in Russian),  VINITI, Moscow, 1999, 361 -- B99, 1 -- 9;
\label{r:ae}\\
\refstepcounter{r}
\noindent [\ther] \rule{2cm}{.3pt}, To the question of the variability
 of solutions\\
\hspace*{3cm} of difference equations of Poincar\'e type (in Russian),\\
\hspace*{3cm} VINITI, Moscow, 2000, 2416 -- B00, 1 -- 22;\label{r:af}\\
\refstepcounter{r}
\noindent [\ther] \rule{2cm}{.3pt}, On linear forms with coefficients
 in ${\mathbb N}\zeta(1+\mathbb N),$\\
\hspace*{3cm} Max-Plank-Institut f\"ur Mathematik,\\
\hspace*{3cm} Bonn, Preprint Series, 2000, 3, 1 -- 13;
\label{r:ag}\\
\refstepcounter{r}
\noindent [\ther] \rule{2cm}{.3pt},
On the Irrationality of Some Quantyties Containing $\zeta (3)$ (in Russian),\\
\hspace*{3cm} Uspekhi Mat. Nauk, 1979, 34, 3(207), 190;
\label{r:di}\\
\refstepcounter{r}
\noindent [\ther] \rule{2cm}{.3pt}, On the Irrationality of Some Quantities
 Containing $\zeta (3),$\\
\hspace*{3cm} Eleven papers translated from the Russian,\\
\hspace*{3cm} American Mathematical Society, 1988, 140, 45 - 56;
\label{r:ah}\\
\refstepcounter{r}
\noindent [\ther] \rule{2cm}{.3pt}, Linear independence over $\mathbb Q$
 of dilogarithms at rational points\\
\hspace*{3cm} (in Russian), UMN, 37 (1982), 179-180;\\
\hspace*{3cm}english transl. in Russ. Math. surveys 37 (1982), 176-177;
\label{r:ai}\\
\refstepcounter{r}
\noindent [\ther] \rule{2cm}{.3pt}, On a measure of the irrationality
 of dilogarithms at rational points\\
\hspace*{3cm} (in Russian), VINITI, 1984, 4345-84, 1 -- 74;
\label{r:bj}\\
\refstepcounter{r}
\noindent [\ther] \rule{2cm}{.3pt}, To the question of the smallness
of some linear forms\\
\hspace*{3cm} (in Russian), VINITI, 1993, 2413-B93, 1 -- 94;
\label{r:ba}\\
\refstepcounter{r}
\noindent [\ther] \rule{2cm}{.3pt}, About linear forms,\\
\hspace*{3cm} whose coefficients are logarithms\\
\hspace*{3cm} of algebraic numbers  (in Russian),\\
\hspace*{3cm} VINITI, 1995, 135-B95, 1 -- 149;
\label{r:bb}\\
\refstepcounter{r}
\noindent [\ther] \rule{2cm}{.3pt}, About systems of vectors, whose
 coordinates\\
\hspace*{3cm} are linear combinations of logarithms of algebraic numbers\\
\hspace*{3cm} with algebraic coefficients  (in Russian),\\
\hspace*{3cm} VINITI, 1994, 3122-B94, 1 -- 158;
\label{r:bc}\\
\refstepcounter{r}
\noindent [\ther] \rule{2cm}{.3pt}, On the linear forms, whose\\
\hspace*{3cm} coefficients are $\mathbb A$ - linear combinations\\
\hspace*{3cm}  of logarithms of $\mathbb A$ - numbers,\\
\hspace*{3cm} VINITI, 1996,  1617-B96, pp. 1 -- 23.
\label{r:bc1}\\
\refstepcounter{r}
\noindent [\ther] \rule{2cm}{.3pt}, On systems of linear forms, whose\\
\hspace*{3cm} coefficients are $\mathbb A$ - linear combinations\\
\hspace*{3cm}  of logarithms of $\mathbb A$ - numbers,\\
\hspace*{3cm} VINITI, 1996,  2663-B96, pp. 1 -- 18.
\label{r:bc2}\\
\refstepcounter{r}
\noindent [\ther] \rule{2cm}{.3pt}, About linear forms, whose coefficients\\
\hspace*{3cm} are $\mathbb Q$-proportional to the number $\log 2, $
 and the values\\
\hspace*{3cm} of $\zeta (s)$ for integer $s$ (in Russian),\\
\hspace*{3cm} VINITI, 1996, 3258-B96, 1 -- 70;
\label{r:bd}\\
\refstepcounter{r}
\noindent [\ther] \rule{2cm}{.3pt}, The lower estimate for some linear forms,\\
\hspace*{3cm} coefficients of which are proportional to the values\\
\hspace*{3cm} of $\zeta (s)$ for integer $s$ (in Russian),\\
\hspace*{3cm} VINITI, 1997, 3072-B97, 1 -- 77;
\label{r:be}\\
\refstepcounter{r}
\noindent [\ther] \rule{2cm}{.3pt},
 On linear forms with coefficients in
${\mathbb N} \zeta(1 + {\mathbb N}) $\\
\hspace*{3cm} Max-Plank-Institut f\"ur Mathematik, Bonn,\\
\hspace*{3cm} Preprint Series, 2000, 3, 1 -- 13;\label{r:cc0}\\
\refstepcounter{r}
\noindent [\ther] \rule{2cm}{.3pt},
 On linear forms with coefficients
 in ${\mathbb N}\zeta(1+\mathbb N)$\\
\hspace*{3cm} (the detailed version,part 1),
 Max-Plank-Institut f\"ur Mathematik,\\
\hspace*{3cm} Bonn, Preprint Series, 2001, 15, 1 -- 20;
\label{r:bf}\\
\refstepcounter{r}
\noindent [\ther] \rule{2cm}{.3pt}, On linear forms with coefficients
 in ${\mathbb N}\zeta(1+\mathbb N)$\\
\hspace*{3cm} (the detailed version,part 2),
 Max-Plank-Institut f\"ur Mathematik,\\
\hspace*{3cm} Bonn, Preprint Series,2001, 104, 1 -- 36;
\label{r:bg}\\
\refstepcounter{r}
\noindent [\ther] \rule{2cm}{.3pt}, On linear forms with coefficients
in ${\mathbb N}\zeta(1+{\mathbb N})$\\
\hspace*{3cm} (the detailed version,part 3),
 Max-Plank-Institut f\"ur Mathematik,\\
\hspace*{3cm} Bonn, Preprint Series, 2002, 57, 1 -- 33;
\label{r:bh}\\
\refstepcounter{r}
\noindent [\ther] \rule{2cm}{.3pt},
On the rank over ${\mathbb Q}$ of some real matrices (in Russian),\\
\hspace*{3cm} VINITI, 1984, 5736-84; 1 -- 29;\label{r:eg}\\
\refstepcounter{r}
\noindent [\ther] \rule{2cm}{.3pt},
On the rank over ${{\mathbb Q}}$ of some real matrices,\\
\hspace*{3cm} Max-Plank-Institut f\"ur Mathematik,\\
\hspace*{3cm} Bonn, Preprint Series, 2002, 27, 1 -- 32;\label{r:eh}\\
\refstepcounter{r}
\noindent [\ther] \rule{2cm}{.3pt}, On linear forms with coefficients
 in ${\mathbb N}\zeta(1+{\mathbb N})$\\
\hspace*{3cm} (the detailed version, part 4),
 Max-Plank-Institut f\"ur Mathematik,\\
\hspace*{3cm} Bonn, Preprint Series, 2002, 142, 1 -- 27;
\label{r:bi}\\
\refstepcounter{r}
\noindent [\ther] \rule{2cm}{.3pt}, On the dimension of some linear spaces\\
\hspace*{3cm} over finite extension of ${\mathbb Q}$ (part 2),\\
\hspace*{3cm} Max-Plank-Institut f\"ur Mathematik, Bonn, Preprint Series,\\
\hspace*{3cm} 2002, 107, 1 -- 37;
\label{r:cj}\\
\refstepcounter{r}
\noindent [\ther] \rule{2cm}{.3pt}, On the dimension of some linear spaces
 over $\mathbb Q$ (part 3),\\
\hspace*{3cm} Max-Plank-Institut f\"ur Mathematik, Bonn,\\
\hspace*{3cm} Preprint Series, 2003, 16, 1 -- 45.
\label{r:ca}\\
\refstepcounter{r}
\noindent [\ther] \rule{2cm}{.3pt},On the difference equation
 of Poincar\'e type (Part 1).\\
\hspace*{3cm} Max-Plank-Institut f\"ur Mathematik, Bonn,\\
\hspace*{3cm} Preprint Series, 2003, 52, 1 -- 44.
\label{r:cb}\\
\refstepcounter{r}
\noindent [\ther] \rule{2cm}{.3pt},  On the dimension of some linear
 spaces over $\mathbb Q,$ (part 4)\\
\hspace*{3cm} Max-Plank-Institut f\"ur Mathematik, Bonn,\\
\hspace*{3cm} Preprint Series, 2003, 73, 1 -- 38.
\label{r:cc}\\
\refstepcounter{r}
\noindent [\ther] \rule{2cm}{.3pt}, On linear forms with coefficients
 in ${\mathbb N}\zeta(1+\mathbb N)$\\
\hspace*{3cm} (the detailed version, part 5),\\
\hspace*{3cm} Max-Plank-Institut f\"ur Mathematik, Bonn,\\
\hspace*{3cm} Preprint Series, 2003, 83, 1 -- 13.
\label{r:cd0}\\
\refstepcounter{r}
\noindent [\ther] \rule{2cm}{.3pt}, On linear forms with coefficients
 in ${\mathbb N}\zeta(1+\mathbb N)$\\
\hspace*{3cm} (the detailed version, part 6),\\
\hspace*{3cm} Max-Plank-Institut f\"ur Mathematik, Bonn,\\
\hspace*{3cm} Preprint Series, 2003, 99, 1 -- 33.
\label{r:ce0}\\
\refstepcounter{r}
\noindent [\ther] \rule{2cm}{.3pt},On the difference equation
of Poincar\'e type (Part 2).\\
\hspace*{3cm} Max-Plank-Institut f\"ur Mathematik, Bonn,\\
\hspace*{3cm} Preprint Series, 2003, 107, 1 -- 25.
\label{r:cf0}\\
\refstepcounter{r}
\noindent [\ther] \rule{2cm}{.3pt},
On the asymptotic behavior of solutions\\
\hspace*{3cm}  of difference equation (in English).\\
\hspace*{3cm} Chebyshevskij sbornik,
 2003, v.4, issue 2, 142 -- 153.
\label{r:cg0}\\
\refstepcounter{r}
\noindent [\ther] \rule{2cm}{.3pt}, On linear forms with coefficients
 in ${\mathbb N}\zeta(1+\mathbb N),$\\
\hspace*{3cm} Bonner Mathematishe Schriften Nr. 360,\\
\hspace*{3cm} Bonn, 2003, 360.
\label{r:ch0}\\
\refstepcounter{r}
\noindent [\ther] \rule{2cm}{.3pt},  On the dimension of some linear
 spaces over $\mathbb Q,$ (part 5)\\
\hspace*{3cm} Max-Plank-Institut f\"ur Mathematik, Bonn,\\
\hspace*{3cm} Preprint Series, 2004, 46, 1 -- 42.
\label{r:ci1}\\
\refstepcounter{r}
\noindent [\ther] \rule{2cm}{.3pt},On the difference equation
of Poincar\'e type (Part 3).\\
\hspace*{3cm} Max-Plank-Institut f\"ur Mathematik, Bonn,\\
\hspace*{3cm} Preprint Series, 2004, 9, 1 -- 33.
\label{r:cf1}\\
\refstepcounter{r}
\noindent [\ther] \rule{2cm}{.3pt}, On linear forms with coefficients
 in ${\mathbb N}\zeta(1+\mathbb N)$\\
\hspace*{3cm} (the detailed version, part 7),\\
\hspace*{3cm} Max-Plank-Institut f\"ur Mathematik, Bonn,\\
\hspace*{3cm} Preprint Series, 2004, 88, 1 -- 27.
\label{r:ci2}
\vskip 10pt
 {\it E-mail:}{\sl\ gutnik$@@$gutnik.mccme.ru}
\end{document}